\documentclass[11pt]{article}%
\usepackage{amsfonts}
\usepackage{amsmath}
\usepackage{geometry}
\usepackage{fancyhdr}
\usepackage{amssymb}
\usepackage{graphicx}
\usepackage{color}%
\providecommand{\U}[1]{\protect\rule{.1in}{.1in}}
\newtheorem{theorem}{Theorem}[section]

\newtheorem{condition}{Condition}[section]

\newtheorem{corollary}[theorem]{Corollary}

\newtheorem{definition}[theorem]{Definition}

\newtheorem{lemma}[theorem]{Lemma}

\newtheorem{proposition}[theorem]{Proposition}

\newenvironment{proof}[1][Proof]{\noindent\textbf{#1.} }{\ \rule{0.5em}{0.5em}}
\begin{document}

\title{Importance Sampling for Multiscale Diffusions}

\author{Paul Dupuis\thanks{Division of Applied Mathematics,
Brown University, Providence, RI 02912 (dupuis@dam.brown.edu). Research supported in part by the National Science Foundation
(DMS-1008331), the Department of Energy (DE-SCOO02413), and the Army Research
Office (W911NF-09-1-0155).} \and Konstantinos Spiliopoulos\thanks{Division of Applied Mathematics,
Brown University, Providence, RI 02912 (kspiliop@dam.brown.edu). Research supported in part
by the Department of Energy (DE-SCOO02413).} \and Hui Wang\thanks{Division of Applied Mathematics,
Brown University, Providence, RI 02912 (huiwang@dam.brown.edu). Research
supported in part by the National Science Foundation (DMS-1008331), the
Department of Energy (DE-SCOO02413).}}

\maketitle

\begin{abstract}
We construct importance sampling schemes for stochastic differential equations
with small noise and fast oscillating coefficients. Standard Monte Carlo
methods perform poorly for these problems in the small noise limit. With
multiscale processes there are additional complications, and indeed the
straightforward adaptation of methods for standard small noise diffusions will
not produce efficient schemes. Using the subsolution approach we construct
schemes and identify conditions under which the schemes will be asymptotically
optimal. Examples and simulation results are provided.
\end{abstract}


\textbf{AMS}: 60F05, 60F10, 60G60

%\pagestyle{myheadings}
%\thispagestyle{plain}
%\markboth{P. DUPUIS, K. SPILIOPOULOS AND H. WANG}{IMPORTANCE SAMPLING FOR MULTISCALE DIFFUSIONS}

\section{Introduction}

\label{S:Introduction}

In this paper we study efficient importance sampling schemes for simulating
rare events associated with the $d$-dimensional stochastic differential
equation (SDE)%
\begin{align}
dX^{\epsilon}(t)  &  =\left[  \frac{\epsilon}{\delta}b\left(  X^{\epsilon
}(t),\frac{X^{\epsilon}(t)}{\delta}\right)  +c\left(  X^{\epsilon}%
(t),\frac{X^{\epsilon}(t)}{\delta}\right)  \right]  dt+\sqrt{\epsilon}%
\sigma\left(  X^{\epsilon}(t),\frac{X^{\epsilon}(t)}{\delta}\right)
dW(t),\nonumber\\
X^{\epsilon}(0)  &  =x, \label{Eq:LDPandA1}%
\end{align}
where $\delta=\delta(\epsilon)\downarrow0$ and
\begin{equation}
\frac{\epsilon}{\delta}\rightarrow\infty\text{ as }\epsilon\downarrow0,
\label{Def:Regime1}%
\end{equation}
and $W(t)$ is a standard $d$-dimensional Wiener process. The functions
$b(x,y),c(x,y)$ and $\sigma(x,y)$ are assumed to be smooth in each variable
and periodic with period $\lambda$ in every direction with respect to the
second variable. The extension to the first order Langevin equation model with
non-periodic random environment will also be discussed.

The need to simulate rare events occurs in many application areas, including
telecommunication, finance, insurance, and computational physics and
chemistry. However, virtually any simulation problem involving rare events
will have a number of mathematical and computational challenges. As it is well
known, standard Monte Carlo sampling techniques perform very poorly in that
the relative errors under a fixed computational effort grow rapidly as the
event becomes more and more rare. Estimating rare event probabilities in the
context of diffusion processes with fast oscillating coefficients presents
extra difficulties due to the additional small parameter $\delta$ and its
interaction with the intensity of the noise $\epsilon$.

A potential application of the methods presented in this paper is to chemical
physics and biology, such as problems involving the folding and binding
kinetics of proteins. These models usually involve rugged potential surfaces
of a complex hierarchical structure with potential minima within potential
minima, separated by barriers of varying heights due to the presence of
multiple energy scales. In some cases, one can approximate the dynamics by a
diffusion in a rough potential where a smooth potential function is
superimposed by a rough function (see Figures \ref{F:Figure1a} and
\ref{F:FigureWithRandomPotential}). A representative, but by no means
complete, list of references is
\cite{Ansari,BryngelsonOnuchicWolynes,HyeonThirumalai,
LifsonJackson,MondalGhosh,SavenWangWolynes, Zwanzig}. It turns out that these
models often can be approximated by homogenized systems where the effect of
the multiscale nature is partially captured by the effective diffusivity of
the system. The formulas for the effective diffusivity in the aforementioned
chemistry and biology literature coincide with those produced by the
approximation via homogenization, which justifies our assumption that
$\epsilon$ and $\delta$ are related according to (\ref{Def:Regime1}). Note
that the condition (\ref{Def:Regime1}) corresponds to Regime $1$ in
\cite{DupuisSpiliopoulos}, where the sample path large deviation properties of
multiscale diffusions are studied under various regimes.

The aim of this paper is to present a more efficient approach to the sampling
problem for multiscale diffusions. Using the large deviation and weak
convergence results from \cite{DupuisSpiliopoulos}, we show how to construct
asymptotically optimal importance sampling schemes with rigorous bounds on
performance. The construction is based on subsolutions for an associated
partial differential equation as in \cite{DupuisWang2}. However, it becomes
applicable only after significant modifications that take into consideration
the multiscale aspect of the model. More precisely, changes of measure that
are purely based on the homogenized system and directly suggested by its
associated partial differential equation do not lead to efficient importance
sampling schemes. Instead, appropriate modifications involving the solution to
a so-called auxiliary \textquotedblleft cell problem\textquotedblright\ have
to be made in order to achieve asymptotic optimality. This is consistent with
the large deviations results obtained in \cite{DupuisSpiliopoulos}, where a
change of measure (or equivalently a control) in partial feedback form has to
be used to prove a large deviation lower bound. By \textquotedblleft partial
feedback\textquotedblright\ we mean that the change of measure is a function
of the fast variable $X^{\epsilon}/\delta$. In the present paper a control in
full feedback form, i.e., a function of both the slow variable $X^{\epsilon}$
and the fast variable $X^{\epsilon}/\delta$, will be used to construct dynamic
importance sampling schemes with precise asymptotic performance bounds.

To the best of our knowledge, this is the first work to address the design of
asymptotically optimal importance sampling schemes for multiscale diffusions.
Related importance sampling problems for regular small noise diffusions
without fast oscillations have been recently considered in
\cite{VandenEijndenWeare}, where the schemes are based on the solution to the
corresponding Hamilton-Jacobi-Bellman (HJB) equation, and in
\cite{GuasoniRobertson}. The present work is also related to the theory of
homogenization of HJB equations
\cite{AlvarezBardi2001,BuckdahnIchihara,Evans2,HorieIshii,
KozlovPiatniskii,LionsSouganidis}.

The paper is organized as follows. In Section \ref{S:Preliminaries} we review
the concept of importance sampling and the role of subsolutions to certain HJB
equation for small noise diffusions without multiscale features. In Section
\ref{S:PreliminariesLDP}, we introduce assumptions, notation and review the
large deviations results that we use for (\ref{Eq:LDPandA1}). Furthermore, we
explain why the standard construction of importance sampling schemes based on
the homogenized system fails in the multiscale setting. The main theorem and
its proof are presented in Section \ref{S:MainTheorem}, where the correct
change of measure is identified. In Section \ref{S:ExampleLangevinEquation} we
apply the general results to first order Langevin equations and derive some
useful explicit formulas. Extensions to an equation with random environment
are discussed in Section \ref{SS:RandomCoefficients}. We report simulation
results in Section \ref{S:SimulationResults} for both the periodic and random
cases in one dimension. The computational challenges that one faces when
simulating trajectories of multiscale diffusions are discussed in Appendix
A.

\section{Importance Sampling and Subsolutions}

\label{S:Preliminaries}

In this section we review some known results on importance sampling for small
noise diffusions without the multiscale feature. In particular, we discuss how
subsolutions to a related HJB equation can be used to design and analyze
importance sampling schemes for such systems. The purpose of these discussions
is not only to introduce some basic concepts in importance sampling and
subsolutions, but also to set the stage for discussions on why the standard
procedure is not directly applicable to multiscale diffusion models.

\subsection{Preliminaries on Importance Sampling}

\label{SS:PreliminariesIS} Let $\left\{  X^{\epsilon},\epsilon>0\right\}  $ be
a $d$-dimensional small noise diffusion for which a sample path large
deviation principle holds, and denote the rate function over the interval
$[t,T]$ by $S_{tT}(\phi)$. Consider a bounded continuous function
$h:\mathbb{R}^{d}\mapsto\mathbb{R}$ and suppose that one is interested in
estimating
\[
\theta(\epsilon)\doteq\mathrm{E}[e^{-\frac{1}{\epsilon}h(X^{\epsilon}%
(T))}|X^{\epsilon}(t)=x]
\]
by Monte Carlo. Define%
\begin{equation}
G(t,x)\doteq\inf_{\phi\in\mathcal{C}([t,T];\mathbb{R}^{d}),\phi(t)=x}\left[
S_{tT}(\phi)+h(\phi(T))\right]  , \label{Eq:VarProb}%
\end{equation}
where $\mathcal{C}([t,T];\mathbb{R}^{d})$ denotes the space of continuous
functions from $[t,T]$ to $\mathbb{R}^{d}$. Then by the contraction principle
\begin{equation}
\lim_{\epsilon\rightarrow0}-\epsilon\log\theta(\varepsilon)=G(t,x).
\label{LDPprinciple1}%
\end{equation}
Let $\Gamma^{\epsilon}(t,x)$ be any unbiased estimator of $\theta(\epsilon)$
that is defined on some probability space with probability measure
$\bar{\mathrm{P}}$. In other words, $\Gamma^{\epsilon}(t,x)$ is a random
variable such that
\[
\bar{\mathrm{E}}\Gamma^{\epsilon}(t,x)=\theta(\epsilon),
\]
where $\bar{\mathrm{E}}$ is the expectation operator associated with
$\bar{\mathrm{P}}$. In this paper we will consider only unbiased estimators.

In Monte Carlo simulation, one generates a number of independent copies of
$\Gamma^{\epsilon}(t,x)$ and the estimate is the sample mean. The specific
number of samples required depends on the desired accuracy, which is measured
by the variance of the sample mean. However, since the samples are independent
it suffices to consider the variance of a single sample. Because of
unbiasedness, minimizing the variance is equivalent to minimizing the second
moment. By Jensen's inequality
\[
\bar{\mathrm{E}}(\Gamma^{\epsilon}(t,x))^{2}\geq(\bar{\mathrm{E}}%
\Gamma^{\epsilon}(t,x))^{2}=\theta(\epsilon)^{2}.
\]
It then follows from (\ref{LDPprinciple1}) that
\[
\limsup_{\epsilon\rightarrow0}-\epsilon\log\bar{\mathrm{E}}(\Gamma^{\epsilon
}(t,x))^{2}\leq2G(t,x),
\]
and thus $2G(t,x)$ is the best possible rate of decay of the second moment.
If
\[
\liminf_{\epsilon\rightarrow0}-\epsilon\log\bar{\mathrm{E}}(\Gamma^{\epsilon
}(t,x))^{2}\geq2G(t,x),
\]
then $\Gamma^{\epsilon}(t,x)$ achieves this best decay rate, and is said to be
\textit{asymptotically optimal}.

We note that even though much of this paper focuses on asymptotically optimal
schemes, asymptotic optimality is not the only practical concern. If optimal
or nearly optimal schemes are too complicated and difficult to implement then
one may prefer to construct non-optimal but simpler schemes. This is also
possible using the subsolution approach that is discussed later in the paper,
and Theorem \ref{T:UniformlyLogEfficientRegime1} identifies a lower bound on
the improvement over ordinary Monte Carlo that will be obtained. In the end,
it is an issue of balance between complexity and feasibility.

\subsection{Large Deviations for Small Noise Diffusions}

\label{SS:LargeDeviationDiffusion} Consider a small noise $d$-dimensional
diffusion process $X^{\epsilon}\doteq\{X^{\epsilon}(s),t\leq s\leq T\}$
satisfying%
\begin{equation}
dX^{\epsilon}(s)=r\left(  X^{\epsilon}(s)\right)  ds+\sqrt{\epsilon}
\Phi\left(  X^{\epsilon}(s)\right)  dW(s),\hspace{0.2cm}X^{\epsilon}(t)=x.
\label{Eq:RegularDiffusion}%
\end{equation}
Throughout this paper we work with the canonical filtered probability space
$(\Omega,\mathfrak{F},\mathrm{P})$ equipped with a filtration $\{\mathfrak{F}
_{t}\}$ that satisfies the usual conditions. Thus $\{\mathfrak{F}_{t}\}$ is
right-continuous and $\mathfrak{F}_{0}$ contains all $\mathrm{P}$-negligible
sets. Since the purpose of this section is expository, we assume for
simplicity that the coefficients $r(x)$ and $\Phi(x)$ are smooth, that the
diffusion matrix
\[
q(x) \doteq\Phi(x) \Phi(x)^{T}%
\]
is uniformly nondegenerate, and that all these functions are uniformly bounded.

We next present a representation theorem proved in \cite{BoueDupuis}, which
will be used here and also later on to analyze importance sampling schemes in
the multiscale setting. Let $\mathcal{A}$ denote the set of all $\mathfrak{F}%
_{t}$-progressively measurable $d$-dimensional processes $v=\left\{
v(s),0\leq s\leq T\right\}  $ that satisfy
\[
\mathrm{E}\int_{0}^{T}\left\Vert v(t)\right\Vert ^{2}dt<\infty.
\]

\begin{theorem}
\label{T:RepresentationTheorem2} Given $\epsilon>0$, let $X^{\epsilon}$ be the
unique strong solution to (\ref{Eq:RegularDiffusion}). Then for any bounded
Borel-measurable function $g$ mapping $\mathcal{C}([t,T];\mathbb{R}^{d})$ into
$\mathbb{R}$,
\[
-\epsilon\log\mathrm{E}\left[  \exp\left\{  -\frac{g(X^{\epsilon})}{\epsilon
}\right\}  \right]  =\inf_{v\in\mathcal{A}}\mathrm{E}\left[  \frac{1}{2}%
\int_{t}^{T}\|v(s)\|^{2}ds+g(X^{\epsilon,v})\right]  ,
\]
where $X^{\epsilon,v}$ is the unique strong solution to the stochastic
differential equation
\begin{equation}
dX^{\epsilon,v}(s)=r\left(  X^{\epsilon,v}(s)\right)  ds+\Phi\left(
X^{\epsilon,v}(s)\right)  \left[  \sqrt{\epsilon}dW(s)+v(s)ds\right]  ,~~~t\le
s\le T, \label{Eq:ControlledSDe}%
\end{equation}
with initial condition $X^{\epsilon,v}(t)=x$.
\end{theorem}

\noindent It is well known that under these conditions the sample path large
deviation principle holds for $\{X^{\epsilon},\epsilon>0\}$ with rate function%
\[
S_{tT}(\phi)=%
\begin{cases}
\displaystyle{\frac{1}{2}\int_{t}^{T}\left\|  \dot{\phi}(s)-r(\phi
(s))\right\|  ^{2}_{q^{-1}(\phi(s))} ds } & \text{if }\phi\in\mathcal{AC}%
([t,T];\mathbb{R}^{d}),\phi(t)=x\\
+\infty & \text{otherwise},
\end{cases}
\]
where $\mathcal{AC}([t,T];\mathbb{R}^{d})$ denotes the collection of
$\mathbb{R}^{d}$-valued absolutely continuous functions on interval $[t,T]$
and
\[
\|v\|_{B} \doteq\sqrt{v^{T} B v}%
\]
for any $v\in\mathbb{R}^{d}$ and symmetric positive definite matrix $B$. When
$B$ is the identity matrix, $\|v\|_{B}$ is just the standard Euclidean norm
$\|v\|$.

\subsection{Importance Sampling in the Absence of Multiscale Features}

\label{SS:RoleOfSubsolutions}We first recall the notion of a subsolution to an
HJB equation of the type
\begin{equation}
U_{t}(t,x)+\bar{H}(x,\nabla_{x}U(t,x))=0,\quad U(T,x)=h(x).
\label{Eq:HJBequation2}%
\end{equation}
In this paper we consider mostly classical sense subsolutions. In some
circumstances other types, such as weak sense subsolutions, may be useful
\cite{DeanDupuis, DupuisWang2}.

\begin{definition}
\label{Def:ClassicalSubsolution} A function $\bar{U}(t,x):[0,T]\times
\mathbb{R}^{d}\mapsto\mathbb{R}$ is a classical subsolution to the HJB
equation (\ref{Eq:HJBequation2}) if

\begin{enumerate}
\item $\bar{U}$ is continuously differentiable,

\item $\bar{U}_{t}(t,x)+\bar{H}(x,\nabla_{x}\bar{U}(t,x))\geq0$ for every
$(t,x)\in(0,T)\times\mathbb{R}^{d}$,

\item $\bar{U}(T,x)\leq h(x)$ for $x\in\mathbb{R}^{d}$.
\end{enumerate}
\end{definition}

When using subsolutions for importance sampling it is often necessary to
impose stronger regularity conditions somewhat beyond those of Definition
\ref{Def:ClassicalSubsolution}. To ease exposition, we will assume the
following condition throughout the paper. It is by no means most economical.
In particular, the uniform bound on the first and second derivatives is not
necessary, and can be replaced by milder conditions with further effort.
However, it is convenient for the purpose of illustration since it guarantees
the feedback control used in importance sampling is uniformly bounded and thus
circumvents a number of technicalities.

\begin{condition}
\label{Cond:ExtraReg} $\bar{U}$ has continuous derivatives up to order $1$ in
$t$ and order $2$ in $x$, and the first and second derivatives in $x$ are
uniformly bounded.
\end{condition}

Next we review the connection between subsolutions and the performance of
related importance sampling schemes. Typically one designs a subsolution for a
specific starting time and initial state $(t,x)$. With an abuse of notation
$(t,x)$ will also be used at times to denote a generic point in $[0,T]\times
\mathbb{R}^{d}$ (the intended use will be clear from the context). The form of
the Hamiltonian is naturally suggested by the calculus of variation problem
(\ref{Eq:VarProb}) and the explicit formula of the rate function $S_{tT}%
(\phi)$ in Section \ref{SS:LargeDeviationDiffusion}:
\begin{equation}
\bar{H}(x,p)=\inf_{u\in\mathbb{R}^{d}}\left[  \left\langle p,r(x)+\Phi
(x)u\right\rangle +\frac{1}{2}\left\Vert u\right\Vert ^{2}\right]
=\left\langle r(x),p\right\rangle -\frac{1}{2}\langle p,q(x)p\rangle.
\label{Eq:ControlFormHJB}%
\end{equation}
In fact, under mild conditions $G$ is the unique viscosity solution to
(\ref{Eq:HJBequation2}). Let $\bar{U}(t,x)$ be a classical subsolution to
(\ref{Eq:HJBequation2}) and $\bar{u}$ the feedback control defined by the
minimizer in (\ref{Eq:ControlFormHJB}) with $p$ replaced by $\nabla_{x}\bar
{U}(t,x)$, i.e.,
\begin{equation}
\bar{u}(t,x)=-\Phi(x)^{T}\nabla_{x}\bar{U}(t,x). \label{Eq:feedback_control}%
\end{equation}
Note that under the given conditions $\bar{u}(t,x)$ is Lipschitz continuous in
$x$, continuous in $(t,x)$, and uniformly bounded.

Consider the family of probability measures $\bar{\mathrm{P}}^{\epsilon}%
$\ defined by the change of measure
\[
\frac{d\bar{\mathrm{P}}^{\epsilon}}{d\mathrm{P}}=\exp\left\{  -\frac
{1}{2\epsilon}\int_{t}^{T}\left\Vert \bar{u}(s,X^{\epsilon}(s))\right\Vert
^{2}ds+\frac{1}{\sqrt{\epsilon}}\int_{t}^{T}\left\langle \bar{u}%
(s,X^{\epsilon}(s)),dW(s)\right\rangle \right\}  .
\]
By Girsanov's Theorem
\[
\bar{W}(s)=W(s)-\frac{1}{\sqrt{\epsilon}}\int_{t}^{s}\bar{u}(\rho,X^{\epsilon
}(\rho))d\rho,~~~t\leq s\leq T
\]
is a Brownian motion on $[t,T]$ under the probability measure $\bar
{\mathrm{P}}^{\epsilon}$, and $X^{\epsilon}$ satisfies $X^{\epsilon}(t)=x$
and
\[
dX^{\epsilon}(s)=r\left(  X^{\epsilon}(s)\right)  ds+\Phi\left(  X^{\epsilon
}(s)\right)  \left[  \sqrt{\epsilon}d\bar{W}(s)+\bar{u}(s,X^{\epsilon
}(s))ds\right]  .
\]
Letting
\[
\Gamma^{\epsilon}(t,x)=\exp\left\{  -\frac{1}{\epsilon}h(X^{\epsilon
}(T))\right\}  \frac{d\mathrm{P}}{d\bar{\mathrm{P}}^{\epsilon}}(X^{\epsilon
}),
\]
it follows easily that under $\bar{\mathrm{P}}^{\epsilon}$, $\Gamma^{\epsilon
}(t,x)$ is an unbiased estimator for $\theta(\epsilon)$. The performance of
this estimator is characterized by the decay rate of its second moment
\begin{equation}
Q^{\epsilon}(t,x;\bar{u})\doteq\bar{\mathrm{E}}^{\epsilon}\left[  \exp\left\{
-\frac{2}{\epsilon}h(X^{\epsilon}(T))\right\}  \left(  \frac{d\mathrm{P}%
}{d\bar{\mathrm{P}}^{\epsilon}}(X^{\epsilon})\right)  ^{2}\right]  .
\label{Eq:2ndMoment1}%
\end{equation}

Following \cite{DupuisWang2}, a verification argument can be used to analyze
$Q^{\epsilon}(t,x;\bar{u})$ as $\epsilon\rightarrow0$. To this end, we need an
alternative expression of $Q^{\epsilon}(t,x;\bar{u})$ that allows us to invoke
the representation in Theorem \ref{T:RepresentationTheorem2}. More precisely,
since $\bar{u}(s,x)$ is bounded and continuous, we can define $X^{\epsilon
,-\bar{u}}$ to be the unique strong solution to the equation
\[
dX^{\epsilon,-\bar{u}}(s)=r\left(  X^{\epsilon,-\bar{u}}(s)\right)
ds+\Phi\left(  X^{\epsilon,-\bar{u}}(s)\right)  \left[  \sqrt{\epsilon
}dW(s)-\bar{u}(s,X^{\epsilon,-\bar{u}}(s))ds\right]
\]
on $[t,T]$ with initial condition $X^{\epsilon,-\bar{u}}(t)=x$. Then by Lemma
\ref{L:MeasureChange} (stated later on in generality sufficient for the
multiscale case),
\[
Q^{\epsilon}(t,x;\bar{u})=\mathrm{E}\exp\left\{  -\frac{2}{\epsilon
}h(X^{\epsilon,-\bar{u}}(T))+\frac{1}{\epsilon}\int_{t}^{T}\left\Vert \bar
{u}(s,X^{\epsilon,-\bar{u}}(s))\right\Vert ^{2}ds\right\}  .
\]
Note that since $\bar{u}$ and $h$ are bounded the exponent in the last display
is uniformly bounded. Hence by Theorem \ref{T:RepresentationTheorem2}%
\begin{align}
\lefteqn{-\epsilon\log Q^{\epsilon}(t,x;\bar{u})}\label{Eq:RepforQ}\\
&  =\inf_{v\in\mathcal{A}}\mathrm{E}\left[  \frac{1}{2}\int_{t}^{T}\left\Vert
v(s)\right\Vert ^{2}ds+2h(X^{\epsilon,-\bar{u},v}(T))-\int_{t}^{T}\left\Vert
\bar{u}(s,X^{\epsilon,-\bar{u},v}(s))\right\Vert ^{2}ds\right]  ,\nonumber
\end{align}
where $X^{\epsilon,-\bar{u},v}$ is the unique strong solution to the equation%
\[
dX^{\epsilon,-\bar{u},v}(s)=r\left(  X^{\epsilon,-\bar{u},v}(s)\right)
ds+\Phi\left(  X^{\epsilon,-\bar{u},v}(s)\right)  \left[  \sqrt{\epsilon
}dW(s)-[\bar{u}(s,X^{\epsilon,-\bar{u},v}(s))-v(s)]ds\right]
\]
on $[t,T]$ with initial condition $X^{\epsilon,-\bar{u},v}(t)=x$.

Fix an arbitrary $v\in\mathcal{A}$ and let $\hat X = X^{\epsilon,-\bar{u},v}$.
Since $\bar{U}(t,x)$ is a classical subsolution and $\bar{u}(t,x)$ is the
minimizer in (\ref{Eq:ControlFormHJB}), it follows that%
\[
\bar{U}_{t}(t,x)+\left\langle \nabla_{x}\bar{U}(t,x),r(x)-\Phi(x)\bar
{u}(t,x)\right\rangle \geq\frac{3}{2}\left\Vert \bar{u}(t,x)\right\Vert ^{2}%
\]
for every $(t,x)\in(0,T)\times\mathbb{R}^{d}$. Hence It\^{o}'s formula and
(\ref{Eq:feedback_control}) give%
\begin{align*}
d\bar{U}(s,\hat X(s))  &  \ge\frac{3}{2}\left\Vert \bar{u}(s,\hat
X(s))\right\Vert ^{2}ds +\left\langle \nabla_{x}\bar{U}(s,\hat X(s)),\Phi(\hat
X(s))v(s)\right\rangle ds\\
&  ~~~~~+\sqrt{\epsilon}\left\langle \nabla_{x}\bar{U}(s,\hat X(s)),\Phi( \hat
X(s))dW(s) \right\rangle \\
&  ~~~~~+\frac{\epsilon}{2}\text{tr}\left[  q( \hat X(s)) \nabla_{xx}\bar
{U}(s,\hat X(s))\right]  ds.
\end{align*}
Integrating the last two terms over $[t,T]$ gives a random variable
$R(\epsilon,v)$ that converges in $L^{2}$ to zero as $\epsilon\rightarrow0$,
uniformly in $v\in\mathcal{A}$. Observing that the second term on the
right-hand-side is $-\langle\bar u(s, \hat X(s)), v(s)\rangle$ and using
$\bar{U}(T,x)\leq h(x)$, one obtains
\begin{equation}
\label{Eq:Alowerbound}h(\hat X(T))-\bar{U}(t,x) \geq\int_{t}^{T}\left[
\frac{3}{2}\left\Vert \bar{u}(s,\hat X(s))\right\Vert ^{2}-\left\langle \bar
u(s,\hat X(s)),v(s)\right\rangle \right]  ds+R(\epsilon,v).
\end{equation}
Now we use the last display to bound one of the two $h(X^{\epsilon,-\bar{u}%
,v}(T))$ terms on the right-hand-side of (\ref{Eq:RepforQ}), yielding the
lower bound
\[
\frac{1}{2}\int_{t}^{T}\left\Vert v(s)-\bar{u}(s,\hat X(s))\right\Vert
^{2}ds+h(\hat X(T))+\bar{U}(t,x)+R(\epsilon,v).
\]
Setting $\bar{v}(s)=v(s)-\bar{u}(s,\hat X(s))$, it follows that $\hat
X=X^{\epsilon,\bar{v}}$ with $X^{\epsilon,\bar{v}}$ defined as in
(\ref{Eq:ControlledSDe}). Since $\bar{v}\in\mathcal{A}$, by Theorem
\ref{T:RepresentationTheorem2}%
\[
\mathrm{E}\left[  \frac{1}{2}\int_{t}^{T}\left\Vert \bar{v}(s)\right\Vert
^{2}ds+h(\hat X(T))\right]  \geq-\epsilon\log\mathrm{E}\exp\left\{  -\frac
{1}{\epsilon}h(X^{\epsilon}(T))\right\}  ,
\]
and therefore%
\begin{align*}
\liminf_{\epsilon\rightarrow0}-\epsilon\log Q^{\epsilon}(t,x;\bar{u})  &
\geq\liminf_{\epsilon\rightarrow0}\inf_{v\in\mathcal{A}}\mathrm{E}\left[
\frac{1}{2}\int_{t}^{T}\left\Vert \bar{v}(s)\right\Vert ^{2}ds+h(\hat
X(T))+R(\epsilon,v)\right]  +\bar{U}(t,x)\\
&  \geq\liminf_{\epsilon\rightarrow0}-\epsilon\log\mathrm{E}\exp\left\{
-\frac{1}{\epsilon}h(X^{\epsilon}(T))\right\}  +\bar{U}(t,x)\\
&  = G(t,x)+\bar{U}(t,x).
\end{align*}

Given that $\bar{U}$ is a subsolution, it is automatic that $\bar{U}(t,x)\leq
G(t,x)$. Thus for the scheme to be asymptotically optimal we need $\bar
{U}(t,x)=G(t,x)$ at the starting point $(t,x)$. The subsolution $\bar
{U}(t,x)=0$ corresponds to standard Monte Carlo (i.e., no change of measure),
and we recover the expected decay rate for that case, which is $G(t,x)$. Note
that if one can obtain a bound on $\mathrm{E}\left[  2R(\epsilon,v)\right]  $
that is uniform in $v\in\mathcal{A}$, then \textit{non-asymptotic} bounds on
the variance can also be obtained.

\section{Large Deviation Properties of Multiscale Diffusions}

\label{S:PreliminariesLDP} In this section we introduce assumptions and
notation, and briefly review the large deviations results for multiscale
diffusions \cite{DupuisSpiliopoulos}. We also revisit the subsolution approach
to importance sampling as discussed in the last section, and identify where
the standard construction breaks down if the multiscale feature of the problem
is not incorporated. Throughout this section we assume a periodic environment,
that is, the functions $b(x,y)$, $c(x,y)$, and $\sigma(x,y)$ are periodic with
period $\lambda$ in every direction with respect to the second variable $y$.
The extension to general random environments but with specialized dynamics,
namely first order Langevin equations, is discussed in Section
\ref{SS:RandomCoefficients}.

\subsection{The Large Deviation Principle}

\label{SS:LargeDeviationMultiscel}

We recall that the SDE of interest is
\begin{align}
dX^{\epsilon}(s)  &  =\left[  \frac{\epsilon}{\delta}b\left(  X^{\epsilon
}(s),\frac{X^{\epsilon}(s)}{\delta}\right)  +c\left(  X^{\epsilon}%
(s),\frac{X^{\epsilon}(s)}{\delta}\right)  \right]  dt+\sqrt{\epsilon}%
\sigma\left(  X^{\epsilon}(s),\frac{X^{\epsilon}(s)}{\delta}\right)
dW(s),\nonumber\\
X^{\epsilon}(t)  &  =x. \label{Eq:LDPandA1a}%
\end{align}
The following condition on (\ref{Eq:LDPandA1a}) will be used whenever the
periodic case is discussed.

\begin{condition}
\label{A:Assumption1}

\begin{enumerate}
\item The functions $b(x,y),\sigma(x,y)$ are continuous and globally bounded,
as are their partial derivatives up to order $2$ in $x$ and order $1 $ in $y$.
The function $c(x,y)$ is bounded and Lipschitz continuous.

\item The diffusion matrix $\sigma(x,y)\sigma(x,y)^{T}$ is uniformly nondegenerate.
\end{enumerate}
\end{condition}

The following condition will also be assumed. In the condition, $\mathbb{T}%
^{d}=[0,\lambda]^{d}$ denotes the $d-$dimensional torus.

\begin{condition}
\label{A:Assumption2} Consider the operator defined for smooth $f:\mathbb{T}%
^{d}\rightarrow\mathbb{R}$ by
\[
\mathcal{L}_{x}f(y)=\left\langle b(x,y),\nabla f(y)\right\rangle +\frac{1}%
{2}\text{\emph{tr}}\left[  \sigma(x,y)\sigma(x,y)^{T}\nabla\nabla f(y)\right]
,
\]
together with periodic boundary conditions in $y$. For any fixed $x$, let
$\mu(dy|x)$ be the unique invariant probability measure corresponding to
$\mathcal{L}_{x}$. Then the drift $b$ satisfies the centering condition (cf.
\cite{BLP})%
\[
\int_{\mathbb{T}^{d}}b(x,y)\mu(dy|x)=0.
\]

\end{condition}

Under Condition \ref{A:Assumption2}, for each $\ell\in\{1,\ldots,d\}$ and $x$
there exists a unique function $\chi_{\ell}(x,y)$ that is twice differentiable
and $\lambda-$periodic in every direction in $y$, and which solves
\begin{equation}
\mathcal{L}_{x}\chi_{\ell}(x,y)=b_{\ell}(x,y),\quad\int_{\mathbb{T}^{d}}%
\chi_{\ell}(x,y)\mu(dy|x)=0. \label{Eq:CellProblem2}%
\end{equation}
For a proof see \cite[Theorem 3.3.4]{BLP}. The equation (\ref{Eq:CellProblem2}%
) is known as a \emph{cell problem}. Let
\[
\chi=(\chi_{1},\ldots,\chi_{d}).
\]
As we shall see below, $\chi$ plays a crucial role in the design of
asymptotically efficient importance sampling schemes for multiscale diffusions.

We state here the sample path large deviations principle for the solution of
(\ref{Eq:LDPandA1a}) derived in \cite{DupuisSpiliopoulos}. Large deviations
principles for special cases of (\ref{Eq:LDPandA1a}) can also be found in
\cite{FS, Baldi}.

\begin{theorem}
\label{T:MainTheorem3} Assume Conditions \ref{A:Assumption1} and
\ref{A:Assumption2}, and let $\{X^{\epsilon},\epsilon>0\}$ be the unique
strong solution to (\ref{Eq:LDPandA1a}). Let
\begin{align*}
r(x)  &  =\int_{\mathbb{T}^{d}}\left(  I+\frac{\partial\chi}{\partial
y}\right)  (x,y)c(x,y)\mu(dy|x),\\
q(x)  &  =\int_{\mathbb{T}^{d}}\left(  I+\frac{\partial\chi}{\partial
y}\right)  (x,y)\sigma(x,y)\sigma(x,y)^{T}\left(  I+\frac{\partial\chi
}{\partial y}\right)  (x,y)^{T}\mu(dy|x),
\end{align*}
where $I$ denotes the identity matrix. Then $\{X^{\epsilon},\epsilon>0\}$
satisfies a large deviations principle with rate function
\[
S_{tT}(\phi)=%
\begin{cases}
\displaystyle{\frac{1}{2}\int_{t}^{T}\left\|  \dot{\phi}(s)-r(\phi
(s))\right\|  ^{2}_{q^{-1}(\phi(s))} ds } & \text{if }\phi\in\mathcal{AC}%
([t,T];\mathbb{R}^{d}),\phi(t)=x\\
+\infty & \text{otherwise}.
\end{cases}
\]

\end{theorem}

Comparing with the rate function for small noise diffusions in Section
\ref{SS:LargeDeviationDiffusion}, it is obvious why $r(x)$ and $q(x)$ are
referred to as the \textquotedblleft effective drift\textquotedblright\ and
\textquotedblleft effective diffusivity\textquotedblright\ in the literature.

\subsection{A Naive Use of Subsolutions for Multiscale Diffusions}

In this section we illustrate the failure of the standard construction of
importance sampling schemes by subsolutions as was outlined in Section
\ref{SS:RoleOfSubsolutions}. Even if one uses a subsolution with the maximum
possible value at the starting point, the scheme can be far from optimal if
the multiscale feature is not incorporated.

The large deviation rate function in Theorem \ref{T:MainTheorem3} is identical
to that of a small noise diffusion (\ref{Eq:RegularDiffusion}) with
dispersion matrix $\Phi(x)$ as long as
\begin{equation}
\Phi(x)\Phi(x)^{T}=q(x). \label{eqn:theta}%
\end{equation}
Note that for a given $q(x)$, the choice of $\Phi(x)$ is not unique. However,
the distribution of the solution to (\ref{Eq:RegularDiffusion}) remains the
same no matter which $\Phi(x)$ is used, and so we fix a Lipschitz continuous
diffusion matrix $\Phi(x)$ for which (\ref{eqn:theta}) holds. Due to the form
of the calculus of variation problem in the rate function, the HJB equation
related to the multiscale diffusion model and the Hamiltonian $\bar{H}$ are
exactly the same as in (\ref{Eq:HJBequation2}) and (\ref{Eq:ControlFormHJB}),
respectively. Therefore, given a subsolution $\bar{U}(t,x)$,
(\ref{Eq:feedback_control}) suggests the control
\[
\bar{u}(t,x)=-\Phi(x)^{T}\nabla_{x}\bar{U}(t,x)
\]
which we now blindly apply to the multiscale diffusion process model.

Suppose that one mimics the steps used in Section \ref{SS:RoleOfSubsolutions}
for the new process model. To simplify notation, as before we temporarily
denote $X^{\epsilon,-\bar{u},v}$ by $\hat{X}$. Then in place of
(\ref{Eq:Alowerbound}) one obtains%
\begin{align}
\lefteqn{h(\hat{X}(T))-\bar{U}(t,x)}\nonumber\\
&  \geq\int_{t}^{T}\left[  \frac{3}{2}\left\Vert \bar{u}(s,\hat{X}%
(s))\right\Vert ^{2}+\left\langle \nabla_{x}\bar{U}(s,\hat{X}(s)),\sigma
\left(  \hat{X}(s),\frac{\hat{X}(s)}{\delta}\right)  v(s)\right\rangle
\right]  ds\nonumber\\
&  ~~~+\int_{t}^{T}\left\langle \nabla_{x}\bar{U}(s,\hat{X}(s)),\left[
\frac{\epsilon}{\delta}b+c\right]  \left(  \hat{X}(s),\frac{\hat{X}(s)}%
{\delta}\right)  -r(\hat{X}(s))\right\rangle
ds\label{Eq:NaiveUseOfSubsolution}\\
&  ~~~+\int_{t}^{T}\left\langle \nabla_{x}\bar{U}(s,\hat{X}(s)),\left[
\Phi(\hat{X}(s))-\sigma\left(  \hat{X}(s),\frac{\hat{X}(s)}{\delta}\right)
\right]  \bar{u}(s,\hat{X}(s))\right\rangle ds\nonumber\\
&  ~~~+R(\epsilon,v),\nonumber
\end{align}
where $R(\epsilon,v)\rightarrow0$ in $L^{2}$, uniformly in $v\in\mathcal{A}$.
The second integral term in the right hand side of
(\ref{Eq:NaiveUseOfSubsolution}) involves $\frac{\epsilon}{\delta}b(\hat
{X}(s),\hat{X}(s)/\delta)$. To deal with the fact that $\epsilon
/\delta\uparrow\infty$ as $\epsilon\downarrow0$, we recall the cell problem
(\ref{Eq:CellProblem2}) and define the function $\psi=(\psi_{1},\cdots
,\psi_{d})$, where for each $\ell\in\{1,\dots,d\}$
\[
\psi_{\ell}(t,x,y)=\chi_{\ell}(x,y)\frac{\partial\bar{U}(t,x)}{\partial
x_{\ell}}.
\]
Applying It\^{o}'s formula to $\psi$ and substituting into
(\ref{Eq:NaiveUseOfSubsolution}), we obtain
\begin{align*}
\lefteqn{h(\hat{X}(T))-\bar{U}(t,x)}\\
&  \geq\int_{t}^{T}\left[  \frac{3}{2}\left\Vert \bar{u}(s,\hat{X}%
(s))\right\Vert ^{2}+\left\langle \nabla_{x}\bar{U}(s,\hat{X}(s)),\left(
I+\frac{\partial\chi}{\partial y}\right)  \sigma\left(  \hat{X}(s),\frac
{\hat{X}(s)}{\delta}\right)  v(s)\right\rangle \right]  ds\\
&  ~~~+\int_{t}^{T}\left\langle \nabla_{x}\bar{U}(s,\hat{X}(s)),\left(
I+\frac{\partial\chi}{\partial y}\right)  c\left(  \hat{X}(s),\frac{\hat
{X}(s)}{\delta}\right)  -r(\hat{X}(s))\right\rangle ds\\
&  ~~~+\int_{t}^{T}\left\langle \nabla_{x}\bar{U}(s,\hat{X}(s)),\left[
\Phi(\hat{X}(s))-\left(  I+\frac{\partial\chi}{\partial y}\right)
\sigma\left(  \hat{X}(s),\frac{\hat{X}(s)}{\delta}\right)  \right]  \bar
{u}(s,\hat{X}(s))\right\rangle ds\\
&  ~~~+R(\epsilon,v),
\end{align*}
where again $R(\epsilon,v)\rightarrow0$ in $L^{2}$, uniformly in
$v\in\mathcal{A}$. If this inequality is inserted into the representation
(\ref{Eq:RepforQ}), then it becomes clear that the desired lower bound will
not follow. Under mild conditions, homogenization can be applied as in
\cite[Theorem 2.7]{DupuisSpiliopoulos} implying that the second integral term
vanishes in the limit. However, regarding the third term, one will need the
inequality
\[
\int_{t}^{T}\left\langle \nabla_{x}\bar{U}(s,\hat{X}(s)),\left[  \left(
I+\frac{\partial\chi}{\partial y}\right)  \sigma\left(  \hat{X}(s),\frac
{\hat{X}(s)}{\delta}\right)  -\Phi(\hat{X}(s))\right]  v(s)\right\rangle
ds\geq0
\]
to hold at least approximately for $\epsilon$ small. While the corresponding
bound held trivially in the case without multiscale if $\Phi$ is chosen as
$\sigma$, it will not hold even approximately here regardless of the choice of
$\Phi$. Indeed, homogenization theory \cite[Theorem 2.7]{DupuisSpiliopoulos}
implies that
\[
\mathrm{E}\int_{t}^{T}\left\langle \nabla_{x}\bar{U}(s,\hat{X}(s)),\left(
\tilde{\sigma}\left(  \hat{X}(s),\frac{\hat{X}(s)}{\delta}\right)
-\int_{\mathbb{T}^{d}}\tilde{\sigma}\left(  \hat{X}(s),y\right)  \mu
(dy|\hat{X}(s))\right)  v(s)\right\rangle ^{2}ds\rightarrow0.
\]
where $\tilde{\sigma}=(I+\partial\chi/\partial y)\sigma$. Therefore in order
for the desired inequality to hold one should choose
\[
\Phi(x)=\int_{\mathbb{T}^{d}}\left(  I+\frac{\partial\chi}{\partial y}\right)
\sigma\left(  x,y\right)  \mu(dy|x).
\]
This is not possible since it violates (\ref{eqn:theta}) in general.

From the preceding discussion it is not difficult to see the fundamental
difficulty that the averaged or effective diffusivity is too crude an
approximation for the corresponding control in importance sampling to be
efficient. The form of $q$ in Theorem \ref{T:MainTheorem3} in fact suggests
the correct control, which is
\begin{equation}
\bar{u}(s,x,y)=-\sigma^{T}(x,y)\left(  I+\frac{\partial\chi}{\partial
y}(x,y)\right)  ^{T}\nabla_{x}\bar{U}(s,x), \label{Eq:OptimalControlRegime1}%
\end{equation}
where $y$ will be replaced by the fast motion $X^{\epsilon}(s)/\delta$ in
implementation. A proof of this assertion will be given in the next section.
Not surprisingly, this form of control is consistent with the control used in
the proof of the large deviation lower bound in \cite{DupuisSpiliopoulos}.

\section{Statement and Proof of the Main Result}

\label{S:MainTheorem}

Before stating and proving the main result, we recapitulate the framework and
notation. We are interested in importance sampling estimator for a functional
of the form
\[
\theta(\epsilon)=\mathrm{E}[e^{-\frac{1}{\epsilon}h(X^{\epsilon}%
(T))}|X^{\epsilon}(t)=x]
\]
where $X^{\epsilon}$ satisfies the SDE (\ref{Eq:LDPandA1a}). According to
Theorem \ref{T:MainTheorem3} the relevant HJB is (\ref{Eq:HJBequation2}) with
the Hamiltonian $\bar{H}$ of the form (\ref{Eq:ControlFormHJB}). Furthermore,
if $h$ is bounded and continuous then
\[
G(t,x)\doteq\inf\left[  S_{tT}(\phi)+h(\phi(T))\right]  =\lim_{\epsilon
\rightarrow0}-\epsilon\log\theta(\epsilon),
\]
where the infimum is taken over all $\phi\in\mathcal{C}([t,T];\mathbb{R}^{d})$
such that $\phi(t)=x$.

Let $\bar{U}$ be a subsolution to the HJB equation with the terminal condition
$h$, and define the control $\bar{u}$ by (\ref{Eq:OptimalControlRegime1}).
Letting
\[
u(s)=\bar{u}(s,X^{\epsilon}(s),X^{\epsilon}(s)/\delta),
\]
it follows from Girsanov's Theorem that
\[
dX^{\epsilon}(s)=\left[  \frac{\epsilon}{\delta}b+c\right]  \left(
X^{\epsilon}(s),\frac{X^{\epsilon}(s)}{\delta}\right)  ds+\sigma\left(
X^{\epsilon}(s),\frac{X^{\epsilon}(s)}{\delta}\right)  \left[  \sqrt{\epsilon
}d\bar{W}(s)+u(s)ds\right]  ,
\]
where $\bar{W}(s)$ is a standard Brownian motion under the probability measure
$\bar{\mathrm{P}}^{\epsilon}$ defined by
\begin{equation}
\frac{d\bar{\mathrm{P}}^{\epsilon}}{d\mathrm{P}}=\exp\left\{  -\frac
{1}{2\epsilon}\int_{t}^{T}\left\Vert u(s)\right\Vert ^{2}ds+\frac{1}%
{\sqrt{\epsilon}}\int_{t}^{T}\left\langle u(s),dW(s)\right\rangle \right\}  .
\label{Eq:2ndRNderiv}%
\end{equation}
The performance measure is then given by the decay rate of the second moment
$Q^{\epsilon}(t,x;\bar{u})$ as defined in (\ref{Eq:2ndMoment1}).

\begin{theorem}
\label{T:UniformlyLogEfficientRegime1} Let $\{X^{\epsilon},\epsilon>0\}$ be
the solution to (\ref{Eq:LDPandA1a}). Consider a bounded and continuous
function $h:\mathbb{R}^{d}\mapsto\mathbb{R}$ and assume Conditions
\ref{Cond:ExtraReg}, \ref{A:Assumption1} and \ref{A:Assumption2}. Let $\bar
{U}(s,x)$ be a subsolution of (\ref{Eq:ControlFormHJB}) and define the control
$\bar{u}(s,x,y)$ by (\ref{Eq:OptimalControlRegime1}). Then
\begin{equation}
\liminf_{\epsilon\rightarrow0}-\epsilon\ln Q^{\epsilon}(t,x;\bar{u})\geq
G(t,x)+\bar{U}(t,x). \label{Eq:GoalRegime1Subsolution}%
\end{equation}
%
%\esquare

\end{theorem}

Theorem \ref{T:UniformlyLogEfficientRegime1} does not cover the important case
of estimating probabilities such as $\mathrm{P}[X^{\epsilon}(T)\in
A|X^{\epsilon}(t)=x]$, since in this case the corresponding function $h$ is
neither bounded nor continuous. Recall that a set $A\subset\mathbb{R}^{d}$ is
called \textit{regular} [with respect to $S_{tT}$ and the initial condition
$(t,x)$] if the infimum of $S_{tT}$ over the closure $\bar{A}$ is the same as
the infimum over the interior $A^{o}$. The following result analogous to
Theorem \ref{T:UniformlyLogEfficientRegime1} holds. Its proof uses an argument
very similar to \cite{DupuisWang2} and is thus omitted.

\begin{proposition}
\label{C:UniformlyLogEfficient1} Let $\{X^{\epsilon},\epsilon>0\}$ be the
solution to (\ref{Eq:LDPandA1a}). Assume Conditions \ref{Cond:ExtraReg},
\ref{A:Assumption1} and \ref{A:Assumption2}. Consider a regular set
$A\subset\mathbb{R}^{d}$ and let
\[
h(x)=%
\begin{cases}
0 & \text{if }x\in A\\
+\infty & \text{if }x\notin A.
\end{cases}
\]
Let $\bar{U}(s,x)$ be a subsolution of (\ref{Eq:HJBequation2}) with the
terminal condition $h$ and define the control $\bar{u}(s,x,y)$ by
(\ref{Eq:OptimalControlRegime1}). Then (\ref{Eq:GoalRegime1Subsolution}) holds.
\end{proposition}

The rest of this section is devoted to the proof of Theorem
\ref{T:UniformlyLogEfficientRegime1}. We first establish an alternative
representation for the performance measure $Q^{\epsilon}$ in terms of bounded
functions, which will allow us to invoke Theorem
\ref{T:RepresentationTheorem2}.

\begin{lemma}
\label{L:MeasureChange} Let $\bar{X}^{\epsilon}\doteq\{X^{\epsilon,-\bar{u}%
},t\leq s\leq T\}$ solve $\bar{X}^{\epsilon}(t)=x$ and
\begin{align*}
d\bar{X}^{\epsilon}(s)  &  =\left[  \frac{\epsilon}{\delta}b+c\right]  \left(
\bar{X}^{\epsilon}(s),\frac{\bar{X}^{\epsilon}(s)}{\delta}\right)  ds\\
&  ~~~~~+\sigma\left(  \bar{X}^{\epsilon}(s),\frac{\bar{X}^{\epsilon}%
(s)}{\delta}\right)  \left[  \sqrt{\epsilon}dW(s)-\bar{u}\left(  s,\bar
{X}^{\epsilon}(s),\frac{\bar{X}^{\epsilon}(s)}{\delta}\right)  ds\right]  .
\end{align*}
Then
\begin{align*}
Q^{\epsilon}(t,x;\bar{u})  &  \doteq\bar{\mathrm{E}}^{\epsilon}\left[
\exp\left\{  -\frac{2}{\epsilon}h(X^{\epsilon}(T))\right\}  \left(
\frac{d\mathrm{P}}{d\bar{\mathrm{P}}^{\epsilon}}(X^{\epsilon})\right)
^{2}\right] \\
&  =\mathrm{E}\exp\left\{  -\frac{2}{\epsilon}h(\bar{X}^{\epsilon}%
(T))+\frac{1}{\epsilon}\int_{t}^{T}\left\Vert \bar{u}\left(  s,\bar
{X}^{\epsilon}(s),\bar{X}^{\epsilon}(s)/\delta\right)  \right\Vert
^{2}ds\right\}  .
\end{align*}

\end{lemma}

\begin{proof}
Since $\bar{u}(s,z,z/\delta)$ is uniformly bounded, it follows from Girsanov's
theorem that
\[
\frac{d\mathrm{Q}}{d\mathrm{P}}\doteq\exp\left\{  -\frac{1}{\sqrt{\epsilon}%
}\int_{t}^{T}\left\langle u(s),dW(s)\right\rangle -\frac{1}{2\epsilon}\int
_{t}^{T}\left\Vert u(s)\right\Vert ^{2}ds\right\}
\]
defines a new probability measure $\mathrm{Q}$ under which
\[
\hat{W}(s)=W(s)+\frac{1}{\sqrt{\epsilon}}\int_{t}^{s}u(\rho)d\rho
\]
is a Brownian motion. Therefore $X^{\epsilon}$ under $\mathrm{Q}$ has the same
distribution as $\bar{X}^{\epsilon}$ under $\mathrm{P}$. This implies
\begin{align*}
\lefteqn{\mathrm{E}\exp\left\{  -\frac{2}{\epsilon}h(\bar{X}^{\epsilon
}(T))+\frac{1}{\epsilon}\int_{t}^{T}\left\Vert \bar{u}(s,\bar{X}^{\epsilon
}(s),\bar{X}^{\epsilon}(s)/\delta)\right\Vert ^{2}ds\right\}  }\\
&  =\mathrm{E}^{\mathrm{Q}}\exp\left\{  -\frac{2}{\epsilon}h(X^{\epsilon
}(T))+\frac{1}{\epsilon}\int_{t}^{T}\Vert u(s)\Vert^{2}ds\right\} \\
&  =\mathrm{E}\exp\left\{  -\frac{2}{\epsilon}h(X^{\epsilon}(T))+\frac
{1}{2\epsilon}\int_{t}^{T}\Vert u(s)\Vert^{2}ds-\frac{1}{\sqrt{\epsilon}}%
\int_{t}^{T}\left\langle u(s),dW(s)\right\rangle \right\}  .
\end{align*}
Using (\ref{Eq:2ndRNderiv}), we can continue the last display as
\begin{align*}
\lefteqn{\bar{\mathrm{E}}^{\epsilon}\exp\left\{  -\frac{2}{\epsilon
}h(X^{\epsilon}(T))+\frac{1}{\epsilon}\int_{t}^{T}\Vert u(s)\Vert^{2}%
ds-\frac{2}{\sqrt{\epsilon}}\int_{t}^{T}\left\langle u(s),dW(s)\right\rangle
\right\}  }\\
&  =\bar{\mathrm{E}}^{\epsilon}\left[  \exp\left\{  -\frac{2}{\epsilon
}h(X^{\epsilon}(T))\right\}  \left(  \frac{d\mathrm{P}}{d\bar{\mathrm{P}%
}^{\epsilon}}(X^{\epsilon})\right)  ^{2}\right]  .\rule{4cm}{0cm}%
\end{align*}
This completes the proof. \hfill
\end{proof}

\vspace{0.5cm}

\begin{proof}
[Proof of Theorem \ref{T:UniformlyLogEfficientRegime1}]Since $h$ and $\bar{u}$
are bounded, it follows from Lemma \ref{L:MeasureChange} and Theorem
\ref{T:RepresentationTheorem2} that
\begin{align}
\lefteqn{-\epsilon\log Q^{\epsilon}(t,x;\bar{u})}\label{Eq:ToBeBounded}\\
&  =\inf_{v\in\mathcal{A}}\mathrm{E}\left[  \frac{1}{2}\int_{t}^{T}\left\Vert
v(s)\right\Vert ^{2}ds-\int_{t}^{T}\Vert\bar{u}(s,\hat{X}(s),\hat{X}%
(s)/\delta)\Vert^{2}ds+2h(\hat{X}(T))\right]  ,\nonumber
\end{align}
where $\hat{X}=\bar{X}^{\epsilon,v}$ solves $\hat{X}(t)=x$ and
\begin{align*}
d\hat{X}(s)  &  =\left[  \frac{\epsilon}{\delta}b\left(  \hat{X}(s),\frac
{\hat{X}(s)}{\delta}\right)  +\bar{c}\left(  s,\hat{X}(s),\frac{\hat{X}%
(s)}{\delta}\right)  \right]  ds\\
&  ~~~+\sigma\left(  \hat{X}(s),\frac{\hat{X}(s)}{\delta}\right)  \left[
\sqrt{\epsilon}dW(s)+v(s)ds\right]  ,
\end{align*}
with%
\[
\bar{c}\left(  s,x,y\right)  =c(x,y)-\sigma(x,y)\bar{u}(s,x,y).
\]

The asymptotic analysis of variational problems analogous to
(\ref{Eq:ToBeBounded}) has already appeared in \cite{DupuisSpiliopoulos},
where large deviation properties of multiscale diffusions such as Theorem
\ref{T:MainTheorem3} have been established through a weak convergence
approach. To be more precise, the condition $\epsilon/\delta\rightarrow\infty$
corresponds to what is called \textit{Regime 1} in \cite{DupuisSpiliopoulos},
and occupation measure techniques are used to characterize the limit of
variational problems. Tightness and characterization in terms of
\textquotedblleft relaxed controls\textquotedblright\ are proved in
Proposition 3.1 and Theorem 2.8 of \cite{DupuisSpiliopoulos}, and then the
relaxed control formulation is rewritten in terms of an ordinary control in
Theorem 5.2.

The difference between the current variational problem and those considered in
\cite{DupuisSpiliopoulos} is that in \cite{DupuisSpiliopoulos} $\bar{c}$ was
independent of time and the middle term in the right-hand-side of
(\ref{Eq:ToBeBounded}) was absent. Nonetheless, these differences are only
superficial and the analysis of quantities similar to the middle term of
(\ref{Eq:ToBeBounded}), e.g.,
\[
\int_{t}^{T}f\left(  s,\hat{X}(s),\frac{\hat{X}(s)}{\delta}\right)  \,ds,
\]
can be carried out using the same arguments as in \cite[Proposition 3.1 and
Theorem 2.8]{DupuisSpiliopoulos}. For this reason, we will directly state a
bound for (\ref{Eq:ToBeBounded}) without giving the details of the analysis:
\begin{align}
\lefteqn{\liminf_{\epsilon\rightarrow0}-\epsilon\log Q^{\epsilon}(t,x;\bar
{u})}\nonumber\\
&  \geq\inf_{\phi\in\mathcal{AC}([t,T];\mathbb{R}^{d}),\phi(t)=x}\left[
\frac{1}{2}\int_{t}^{T}\left\Vert \dot{\phi}(s)-\bar{r}(s,\phi(s))\right\Vert
_{q^{-1}(\phi(s))}^{2}ds\right. \label{Eq:VarLim}\\
&  \hspace{3cm}\left.  -\int_{t}^{T}\int_{\mathbb{T}^{d}}\left\Vert \bar
{u}(s,\phi(s),y)\right\Vert ^{2}\mu(dy|\phi(s))ds+2h(\phi(T))\right]
,\nonumber
\end{align}
where
\[
\bar{r}(s,x)\doteq r(x)-\int_{\mathbb{T}^{d}}\left(  I+\frac{\partial\chi
}{\partial y}(x,y)\right)  \sigma(x,y)\bar{u}(s,x,y)\mu(dy|x).
\]
Recalling the definition of $\bar{u}$,%
\[
\bar{r}(s,x)=r(x)+q(x)\nabla_{x}\bar{U}(s,x)
\]
and
\[
\int_{t}^{T}\int_{\mathbb{T}^{d}}\left\Vert \bar{u}(s,\phi(s),y)\right\Vert
^{2}\mu(dy|\phi(s))ds=\int_{t}^{T}\langle\nabla_{x}\bar{U}(s,\phi
(s)),q(x)\nabla_{x}\bar{U}(s,\phi(s))\rangle ds.
\]
Thus the quantity to be minimized in (\ref{Eq:VarLim}) can be rewritten as
\begin{align}
\lefteqn{\frac{1}{2}\int_{t}^{T}\left\Vert \dot{\phi}(s)-r(\phi(s))\right\Vert
_{q^{-1}(\phi(s))}^{2}ds-\int_{t}^{T}\langle\dot{\phi}(s)-r(\phi
(s)),\nabla_{x}\bar{U}(s,\phi(s))\rangle ds}\label{Eq:VarLim2}\\
&  ~~~~~-\frac{1}{2}\int_{t}^{T}\langle\nabla_{x}\bar{U}(s,\phi(s)),q(x)\nabla
_{x}\bar{U}(s,\phi(s))\rangle ds+2h(\phi(T)).\rule{2cm}{0cm}\nonumber
\end{align}
Given an arbitrary $\phi\in\mathcal{AC}([t,T];\mathbb{R}^{d})$ with
$\phi(t)=x$, the subsolution property implies that
\begin{align*}
\frac{d}{ds}\bar{U}(s,\phi(s))  &  =\bar{U}_{t}(s,\phi(s))+\langle\nabla
_{x}\bar{U}(s,\phi(s)),\dot{\phi}(s)\rangle\\
&  \geq\langle\nabla_{x}\bar{U}(s,\phi(s)),\dot{\phi}(s)-r(\phi(s))\rangle
+\frac{1}{2}\langle\nabla_{x}\bar{U}(s,\phi(s)),q(\phi(s))\nabla_{x}\bar
{U}(s,\phi(s))\rangle.
\end{align*}
Integrating both sides on $[t,T]$ and using the terminal condition $\bar
{U}(T,x)\leq h(x)$, it follows that (\ref{Eq:VarLim2}) is bounded from below
by
\[
\frac{1}{2}\int_{t}^{T}\left\Vert \dot{\phi}(s)-r(\phi(s))\right\Vert
_{q^{-1}(\phi(s))}^{2}ds+h(\phi(T))+\bar{U}(t,x).
\]
Note that the first summand is $S_{tT}(\phi)$ and by definition $G(t,x)$ is
the infimum of the sum of the first two terms over $\phi\in\mathcal{AC}%
([t,T];\mathbb{R}^{d})$ with $\phi(t)=x$, and thus%
\[
\liminf_{\epsilon\rightarrow0}-\epsilon\log Q^{\epsilon}(t,x;\bar{u})\geq
G(t,x)+\bar{U}(t,x).
\]
This concludes the proof. \hfill
\end{proof}

\section{First Order Langevin Equation with Periodic Environment}

\label{S:ExampleLangevinEquation}

\begin{figure}
[ptb]
\begin{center}
\includegraphics[
height=2.8253in,
width=3.1782in
]%
{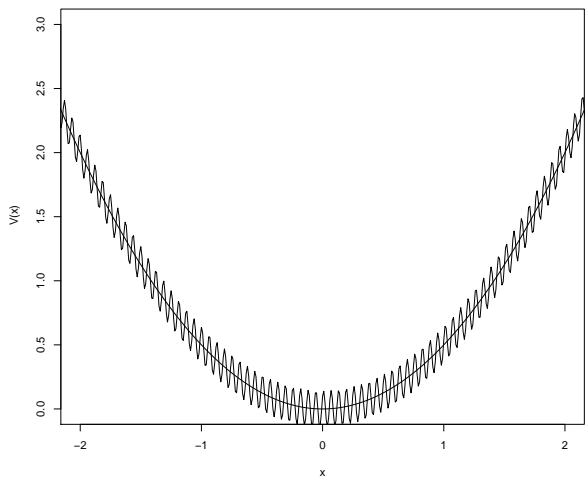}%
\caption{$V^{\epsilon}(x,x/\delta)=\epsilon(\cos(x/\delta)+\sin(x/\delta
))+x^{2}/2$ with $\epsilon=0.1,\delta=0.01$.}%
\label{F:Figure1a}%
\end{center}
\end{figure}

In this section, we apply the general results to a special but important class
of diffusion models, namely, the first order Langevin equation
\[
dX^{\epsilon}(t)=-\nabla V^{\epsilon}\left(  X^{\epsilon}(t),\frac
{X^{\epsilon}(t)}{\delta}\right)  dt+\sqrt{\epsilon}\sqrt{2D}dW(t),\hspace
{0.2cm}X^{\epsilon}(0)=x,
\]
where $V^{\epsilon}$ is some potential function and $2D$ the diffusion
constant. We are particularly interested in the case where the potential
function $V^{\epsilon}$ is composed of a large-scale smooth part and a fast
oscillating part of smaller magnitude:
\[
V^{\epsilon}\left(  x,x/\delta\right)  =\epsilon Q(x/\delta)+V(x).
\]
Thus the equation of interest is%
\begin{equation}
dX^{\epsilon}(t)=\left[  -\frac{\epsilon}{\delta}\nabla Q\left(
\frac{X^{\epsilon}(t)}{\delta}\right)  -\nabla V\left(  X^{\epsilon
}(t)\right)  \right]  dt+\sqrt{\epsilon}\sqrt{2D}dW(t),\hspace{0.2cm}%
X^{\epsilon}(0)=x. \label{Eq:LangevinEquation2}%
\end{equation}
In the notation of previous sections, this corresponds to
\[
b(x,y)=-\nabla Q(y),~~~c(x,y)=-\nabla V(x),~~~\sigma(x,y)=\sqrt{2D}.
\]
An example of such a potential is given in Figure \ref{F:Figure1a}. As before,
we examine in some detail the model (\ref{Eq:LangevinEquation2}) with periodic
environment, and it is assumed in this section that $Q(y)$ is periodic with
period $\lambda$. This periodicity assumption may seem too artificial in many
practical applications. However, it motivates by analogy the design of
importance sampling schemes for first order Langevin equations with general
random environment. See Section \ref{SS:RandomCoefficients} for a discussion
on these extensions.%

%EndExpansion

An important observation for an equation of the form
(\ref{Eq:LangevinEquation2}) is that the cell problem (\ref{Eq:CellProblem2})
depends only on $y$ and not on $x$. Hence, in order to compute $\bar{u}$ from
(\ref{Eq:OptimalControlRegime1}) we need only solve the cell problem
(\ref{Eq:CellProblem2}) \textit{once}. To be more specific, the invariant
distribution $\mu(dy|x)$ to the cell problem is independent of $x$, is of
Gibbs type
\[
\mu(dy|x)=\mu(dy) = \frac{1}{L}e^{-\frac{Q(y)}{D}}dy,\hspace{0.2cm}%
~~L=\int_{\mathbb{T}^{d}}e^{-\frac{Q(y)}{D}}dy,
\]
and satisfies Condition \ref{A:Assumption2}.

Explicit formulas for the large deviation rate functions and related
quantities are readily available for the one dimensional case $d=1$. For
multi-dimensional cases, they are also available under extra assumptions on
the potential function $V$ \cite{DupuisSpiliopoulos}. Since our numerical
simulation will be performed on one-dimensional models, we only state the
relevant results in Corollary \ref{C:MainCorollary1} and refer the readers to
\cite{DupuisSpiliopoulos} for more general formulas. The proof is omitted as
it is a straightforward calculation from Theorems \ref{T:MainTheorem3} and
\ref{T:UniformlyLogEfficientRegime1}.

\begin{corollary}
\label{C:MainCorollary1} Consider the one dimensional case $d=1$ and let
$\{X^{\epsilon}\}$ be the unique strong solution to
(\ref{Eq:LangevinEquation2}). Under Condition \ref{A:Assumption1},
$\{X^{\epsilon}\}$ satisfies a large deviations principle with rate function
\begin{equation}
S_{0T}(\phi)=%
\begin{cases}
\displaystyle{\frac{1}{2}\int_{0}^{T}\frac{1}{q}[\dot{\phi}(s)-r(\phi
(s))]^{2}ds} & \text{if }\phi\in\mathcal{AC}([0,T];\mathbb{R})\text{ and }%
\phi(0)=x\\
+\infty & \text{otherwise},
\end{cases}
\label{ActionFunctional1_1}%
\end{equation}
where
\[
r(x)=-\frac{\lambda^{2}V^{\prime}(x)}{L\hat{L}},\hspace{0.5cm}q=\frac
{2D\lambda^{2}}{L\hat{L}}%
\]
and
\[
L=\int_{\mathbb{T}}e^{-\frac{Q(y)}{D}}dy,\hspace{0.5cm}\hat{L}=\int
_{\mathbb{T}}e^{\frac{Q(y)}{D}}dy.
\]
Given a classical subsolution $\bar{U}$, the importance sampling control that
appears in Theorem \ref{T:UniformlyLogEfficientRegime1} takes the form
\begin{equation}
\bar{u}(t,x,y)=-\frac{\sqrt{2D}\lambda}{\hat{L}}e^{\frac{Q(y)}{D}}\bar{U}%
_{x}(t,x). \label{Eq:OptimalControlExample1}%
\end{equation}

\end{corollary}

An interesting observation from Corollary \ref{C:MainCorollary1} is that the
effective diffusivity $\epsilon q$ is always smaller than the diffusivity
$\epsilon2D$ of the unhomogenized equation, since by H\"{o}lder's inequality
\[
L\hat{L}>\left(  \int_{\mathbb{T}}dy\right)  ^{2}=\lambda^{2},
\]
as long as $Q$ is not a constant. The intuition is that the potential surface
has many small local minima, which manifest themselves in the homogenized
dynamics by a reduction in the diffusion coefficient since a particle
traveling on a rough potential surface may suffer from the \textquotedblleft
trapping\textquotedblright\ effect of these local minima.

\section{Extension to Random Environment}

\label{SS:RandomCoefficients}

Up until now all the multiscale diffusion models we have considered are of
periodic environment, i.e., the drift vector and the dispersion matrix are
both periodic with respect to the fast variable. This section discusses an
extension to a random environment. To illustrate the main idea, we specialize
again to diffusions governed by the first order Langevin equation of type
\begin{equation}
dX^{\epsilon,\delta}(t)=\left[  -\frac{\epsilon}{\delta}\nabla Q\left(
\frac{X^{\epsilon}(t)}{\delta}\right)  -\nabla V\left(  X^{\epsilon
}(t)\right)  \right]  dt+\sqrt{\epsilon}\sqrt{2D}dW(t),\hspace{0.2cm}%
X^{\epsilon,\delta}(0)=x. \label{Eq:LangevinEquation3}%
\end{equation}
In this section, we find it convenient to keep track of both $\epsilon$ and
$\delta$ and write $X^{\epsilon,\delta}(t)$ for the solution to
(\ref{Eq:LangevinEquation3}) as opposed to $X^{\epsilon}$ as in Section
\ref{S:ExampleLangevinEquation}. The following condition is assumed throughout
this section as a substitute for Condition \ref{A:Assumption1} in the periodic case.

\begin{condition}
\label{A:AssumptionRandomCase}

\begin{enumerate}
\item {The coefficient $\{Q(y), y\in\mathbb{R}^{d}\}$ is a stationary, ergodic
random field defined on some probability space $(\Psi,\mathcal{G},\nu)$. For
every $\omega\in\Psi$, $Q(y,\omega)$ is $\mathcal{C}^{2}(\mathbb{R}^{d})$ in
$y$ with bounded and Lipschitz continuous derivatives up to order 2. }

\item {The coefficient $V(x)$ is deterministic and $V(x)\in\mathcal{C}%
^{2}(\mathbb{R}^{d})$ with bounded and Lipschitz continuous derivatives up to
order 2.}
\end{enumerate}
\end{condition}

\begin{figure}
[ptb]
\begin{center}
\includegraphics[
height=2.7371in,
width=2.744in
]%
{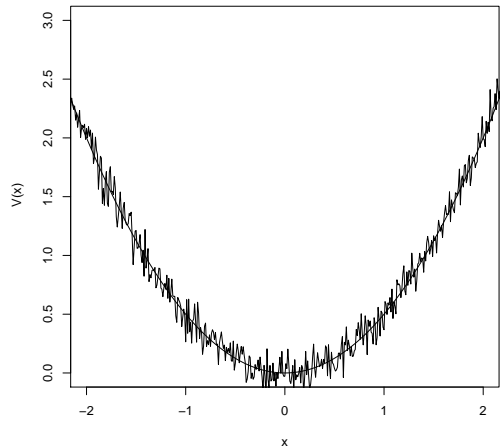}%
\caption{$V^{\epsilon}(x,x/\delta)=\epsilon Q(x/\delta)+x^{2}/2$ with
$\epsilon=0.1$ and $\delta=0.01$.}%
\label{F:FigureWithRandomPotential}%
\end{center}
\end{figure}

The Wiener process in (\ref{Eq:LangevinEquation3}) is defined on another
probability space, and we work with the product space and product measure and
so $W$ is independent of $Q$. Note that for the sake of notational simplicity,
we have suppressed the dependence of $X^{\epsilon,\delta}$ and $Q$ on $\omega$
for $\omega\in\Psi$. Under Condition \ref{A:AssumptionRandomCase}, for every
$\omega\in\Psi$, there exists a unique strong solution to the SDE
(\ref{Eq:LangevinEquation3}). In contrast to the periodic case, when equation
(\ref{Eq:LangevinEquation3}) is used to model the dynamics of a particle in a
rough potential, the roughness is due to the \textquotedblleft
small\textquotedblright\ randomness generated by the random field $\epsilon
Q(x/\delta)$. Figure \ref{F:FigureWithRandomPotential} depicts a realization
of such a random field superimposed on the smooth potential function
$V(x)=x^{2}/2$. In the figure, $Q(x)$ is a zero mean Gaussian random field
with Gaussian type correlation, i.e., $E^{\nu}[Q(x)Q(y)]=\exp[-|x-y|^{2}]$.%

%EndExpansion

Compared with multiscale diffusions with periodic environment, the analysis
for general random media is relatively new. To the best of our knowledge, the
first attempt to generalize the results obtained for the locally periodic
setting to the locally stationary setting probably appeared in the work of
\cite{OllaSiri2004}, which studied random walks on $\mathbb{Z}$ with a locally
stationary environment. Extensions have been considered in
\cite{Rhodes2009a,Rhodes2009b} to diffusions whose generators are self-adjoint
and take a certain form.

\subsection{Homogenization in One Dimension}

\label{S:HomTheorem1D}

As in Section \ref{S:ExampleLangevinEquation}, we state the homogenization
theorem for the one-dimensional case where everything can be explicitly
quantified. For higher dimensions, analogous results exist but explicit
calculation is much more difficult. In the following result we assume that
$\epsilon$ is fixed and let $\delta$ tend to zero in order to clearly identify
the effect of homogenization. To ease exposition, we temporarily denote
$X^{\epsilon,\delta}$ by $X^{\delta}$.

\begin{theorem}
\label{T:HomogenizationThm1D} Consider the one dimension case and assume
Condition \ref{A:AssumptionRandomCase}. Then the law of $X^{\delta}$ on
$\mathcal{C}([0,T];\mathbb{R})$ converges weakly to the law of $X$, in
probability with respect to $\nu$, where $X$ is the solution to the SDE
\[
dX(t)=-\frac{1}{K\hat{K}}V^{\prime}\left(  X(t)\right)  dt+\sqrt
{\frac{2\epsilon D}{K\hat{K}}}dW(t),\hspace{0.2cm}X(0)=x,
\]
and with
\begin{equation}
K=\mathrm{E}^{\nu}\left[  e^{-Q(y)/D}\right]  ,~~~~\hat{K}=\mathrm{E}^{\nu
}\left[  e^{Q(y)/D}\right]  . \label{Eq:HomogenizedCoefficients}%
\end{equation}

\end{theorem}

The rest of this subsection is devoted to the proof of this theorem. Since it
is very similar to that of \cite[Theorem 3.1]{Rhodes2009a}, we will only give
an outline. Without loss of generality, we assume $\epsilon=1$ in the proof.
We start with the following lemma.

\begin{lemma}
\label{L:AveragingRandomMedia} Let $h:\mathbb{R}\rightarrow\mathbb{R}$ be a
bounded and measurable function and let $\{\psi(y),y\in\mathbb{R}\}$ be a
stationary random field such that
\[
E^{\nu}\left[  |\psi(y)|e^{-Q(y)/D}\right]  <\infty.
\]
Then as $\delta\rightarrow0$
\[
E_{\omega}\sup_{0\leq t\leq T}\left\vert \int_{0}^{t}\psi\left(
\frac{X^{\delta}(s)}{\delta}\right)  h(X^{\delta}(s))ds-B\int_{0}%
^{t}h(X^{\delta}(s))ds\right\vert \rightarrow0
\]
in $L^{1}(\nu)$, where $E_{\omega}$ denotes the expected value with respect to the
independent Wiener process $W$ but with $\omega\in\Psi$ given, and
\[
B\doteq\frac{1}{K}E^{\nu}\left[  \psi(y)e^{-Q(y)/D}\right]  .
\]

\end{lemma}

\begin{proof}
The proof follows from Theorem 6.1 in \cite{Rhodes2009a}. The only observation
we need to make is that in our case the fast motion is governed by the
generator
\[
-Q^{\prime}(y) \frac{d}{dy} + D \frac{d^{2}}{dy^{2}}.
\]
Then by the discussion in Section 2.2, in particular Theorem 2.1, of
\cite{Olla1994}, the averaging should be taken under the corresponding ergodic
stationary probability measure on $(\Psi,\mathcal{G})$ say $\pi$ which
satisfies
\[
E^{\pi}[\psi(y)] = \frac{1}{K}E^{\nu}[\psi(y) e^{-Q(y)/D}]
\]
for any stationary random field $\psi$. \hfill
\end{proof}

\vspace{0.4cm}

\begin{proof}
[Proof of Theorem \ref{T:HomogenizationThm1D}]Let $\{\chi(y),y\in\mathbb{R}\}$
be the random field defined by
\[
\chi(y)=\frac{1}{\hat{K}}\int_{0}^{y}e^{Q(z)/D}dz-y
\]
It is straightforward to verify that $\mathrm{E}^{\nu}\left[  \chi(y)\right]
=\mathrm{E}^{\nu}\left[  \chi^{\prime}(y)\right]  =0$ and
\begin{align}
1+\chi^{\prime}(y)  &  =\frac{1}{\hat{K}}e^{Q(y)/D}%
,\label{Eq:AuxillaryProcess}\\
-Q^{\prime}(y)\chi^{\prime}(y)+D\chi^{\prime\prime}(y)  &  =Q^{\prime
}(y).\nonumber
\end{align}

Define the process
\[
Y^{\delta}(t)=X^{\delta}(t)+\int_{0}^{t}\left[  1+\chi^{\prime}\left(
\frac{X^{\delta}(s)}{\delta}\right)  \right]  V^{\prime}(X^{\delta
}(s))ds+\delta\cdot\chi\left(  \frac{X^{\delta}(t)}{\delta}\right)
\]
Then It\^{o}'s formula implies that
\[
Y^{\delta}(t)=x+\sqrt{2D}\int_{0}^{t}\left[  1+\chi^{\prime}\left(
\frac{X^{\delta}(s)}{\delta}\right)  \right]  \,dW(s).
\]
The latter and (\ref{Eq:AuxillaryProcess}) imply that $\{Y^{\delta}\}$ is
$\nu$-a.s. a martingale with quadratic variation
\[
\langle Y^{\delta}\rangle(t)=\frac{2D}{\hat{K}^{2}}\int_{0}^{t}e^{\frac
{2Q(X^{\delta}(s)/\delta)}{D}}ds.
\]
It follows now from Lemma \ref{L:AveragingRandomMedia} [taking $h(x)=1$] that
as $\delta\rightarrow0$
\[
E_{\omega}\sup_{0\leq t\leq T}\left\vert \langle Y^{\delta}\rangle(t)-Bt\right\vert
\rightarrow0\mbox{~~~in $L^1(\nu)$},
\]
where
\[
B\doteq\frac{2D}{\hat{K}^{2}}\cdot\frac{1}{K}\mathrm{E}^{\nu}\left[
e^{Q(y)/D}\right]  =\frac{2D}{\hat{K}^{2}}\cdot\frac{\hat{K}}{K}=\frac
{2D}{K\hat{K}}.
\]

Similarly, by Lemma \ref{L:AveragingRandomMedia} [taking $h(x)=V^{\prime}(x)
$] again we have that
\[
E_{\omega}\sup_{0\leq t\leq T}\left\vert \int_{0}^{t}\left[  1+\chi^{\prime}\left(
\frac{X^{\delta}(s)}{\delta}\right)  \right]  V^{\prime}(X^{\delta
}(s))ds-\frac{1}{K\hat{K}}\int_{0}^{t}V^{\prime}(X^{\delta}(s))ds\right\vert
\rightarrow0
\]
in $L^{1}(\nu)$ as $\delta\rightarrow0$. Finally, note that by Birkhoff's
ergodic theorem $\frac{\chi(y)}{y}\rightarrow0$ $\nu$-almost surely, as
$|y|\rightarrow\infty$. Thus as $\delta\rightarrow0$
\[
\sup_{0\leq t\leq T}\left\vert \delta\cdot\chi\left(  \frac{X^{\delta}%
(t)}{\delta}\right)  \right\vert \rightarrow0\mbox{~~with probability one}
\]
$\nu$-almost surely. Then the desired convergence follows immediately since
\[
Y^{\delta}(\cdot)\Rightarrow x+\sqrt{\frac{2D}{K\hat{K}}}W(\cdot)
\]
in probability with respect to $\nu$. This completes the proof. \hfill
\end{proof}

\subsection{Importance Sampling Schemes}

\label{SS:ISRandomMedia}

It is easy to see that the homogenization result in the random case is
analogous to that in the periodic case Corollary \ref{C:MainCorollary1} with
$\lambda=1$ and $L,\hat{L}$ replaced by $K,\hat{K}$, respectively. By analogy,
it suggests that in the one dimensional case $\{X^{\epsilon,\delta}%
,\epsilon,\delta>0\}$ satisfies the large deviations principle with rate
function given by (\ref{ActionFunctional1_1}), where
\[
r(x)=-\frac{V^{\prime}(x)}{K\hat{K}},\hspace{0.2cm}q=\frac{2D}{K\hat{K}}%
\]
and $K$ and $\hat{K}$ are defined by (\ref{Eq:HomogenizedCoefficients}).
Therefore it is natural to conjecture that given a classical subsolution
$\bar{U}$, the corresponding control takes a similar form
\[
\bar{u}(t,x,y)=-\sqrt{2D}\frac{1}{\hat{K}}e^{\frac{Q(y)}{D}}\bar{U}_{x}(t,x).
\]
Note that in contrast to the periodic case, the control $\bar{u}$ is random in
that it implicitly depends on $\omega\in\Psi$ since the random environment
$Q(y)$ depends on it. Even though the associated importance sampling schemes
are shown to be efficient in our empirical study, development of the
underlying large deviation theory and a rigorous performance analysis remains
to be done.

\section{Simulation Results}

\label{S:SimulationResults} In this section we test the performance of various
Monte Carlo estimators for multiscale diffusions with periodic or random
environment. Throughout this section, we assume that the diffusion process
$X^{\epsilon}$ is one-dimensional and satisfies the first order Langevin
equation (\ref{Eq:LangevinEquation2}).

\subsection{Simulation Results for Periodic Case}

\label{SS:SimulationResultsPeriodicCase} Suppose that we are interested in the
Monte Carlo estimation of
\[
\theta(\epsilon)\doteq\mathrm{E}[e^{-\frac{1}{\epsilon}h(X^{\epsilon}%
(T))}|X^{\epsilon}(t)=x]
\]
for a continuous function $h$ and in a periodic environment. We compare three
unbiased estimators.

\begin{enumerate}
\item The standard Monte Carlo estimator $\hat\theta_{0}(\epsilon)$,

\item The importance sampling estimator $\hat{\theta}_{1}(\epsilon)$ based on
the change of measure suggested by (\ref{Eq:OptimalControlExample1}), i.e.,
\[
u_{1}(t,x,y)\doteq-\sqrt{2D}\frac{\lambda}{\hat{L}}e^{Q(y)/D}\bar{U}%
_{x}(t,x),
\]
where $\bar{U}$ is a subsolution to the homogenized HJB equation
(\ref{Eq:HJBequation2}). Note that in implementation, $x$ will be replaced by
the current state of $X^{\epsilon}$ and $y$ by $X^{\epsilon}/\delta$.

\item The importance sampling estimator $\hat{\theta}_{2}(\epsilon)$ based on
the change of measure corresponding to the control
\[
u_{2}(t,x)\doteq-\sqrt{q}\bar{U}_{x}(t,x)=-\sqrt{2D}\frac{\lambda}{\sqrt
{L\hat{L}}}\bar{U}_{x}(t,x).
\]
This is the change of measure based on the control suggested by the
homogenized HJB equation \textit{without} taking into consideration of the
multiscale nature of the dynamics. It is independent of the fast variable and
differs from the control $u_{1}$ by a factor $\sqrt{2D/q}[1+\chi^{\prime}(y)]$.
\end{enumerate}

Based on the theory developed previously, the estimator $\hat{\theta}%
_{1}(\epsilon)$ should outperform the other two estimators as $\epsilon
\rightarrow0$.

For numerical experimentation, we consider the potential function that is
drawn in Figure \ref{F:Figure1a}, that is, $V(x)=x^{2}/2$ and $Q(y)=\cos
y+\sin y$. Thus the period is $\lambda=2\pi$. Let
\[
h(x)=%
\begin{cases}
(x-1)^{2} & x\geq0\\
(x+1)^{2} & x<0.
\end{cases}
\]
It follows from Corollary \ref{C:MainCorollary1} that the effective drift
$r(x)$ and diffusivity $q$ are
\[
r(x)=-\kappa x\text{,}~~~~q=2D\kappa,~~\mbox{where }\kappa=\frac{\lambda^{2}%
}{L\hat{L}}=\frac{4\pi^{2}}{L\hat{L}}.
\]
The limiting HJB equation (\ref{Eq:HJBequation2}) becomes
\begin{equation}
U_{t}(t,x)-\kappa xU_{x}(t,x)-\kappa D|U_{x}(t,x)|^{2}=0,\quad U(T,x)=h(x).
\label{Eq:1DHJBeqn}%
\end{equation}
Under mild conditions that are satisfied here, the unique viscosity solution
to this HJB equation equals
\[
G(t,x)\doteq\inf\left[  S_{tT}(\phi)+h(\phi(T))\right]  =\lim_{\epsilon
\rightarrow0}-\epsilon\log\theta(\epsilon),
\]
where $S_{tT}(\cdot)$ is given by Corollary \ref{C:MainCorollary1} and the
infimum is taken over all $\phi\in\mathcal{AC}([t,T];\mathbb{R}^{d})$ such
that $\phi(t)=x$. One can solve this variational problem explicitly and
obtain
\[
G(t,x)=\frac{(e^{\kappa T}-|x|e^{\kappa t})^{2}}{(1+2D)e^{2\kappa
T}-2De^{2\kappa t}}.
\]
Since $G$ is not smooth at $x=0$, it is not a classical sense solution. In
general one should mollify it in order to produce a smooth subsolution, but it
is known (see \cite{VandenEijndenWeare} for an analogous situation) that the
bound on performance is still valid if the subsolution is the minimum of two
classical sense solutions with a single discontinuous interface. Therefore we
can just define the subsolution $\bar{U}$ as the solution $G$ and the
corresponding controls are
\begin{align*}
u_{1}(t,x,y)  &  =\sqrt{2D}\frac{2\pi}{\hat{L}}e^{(\cos y+\sin y)/D}%
\frac{2e^{\kappa t}(e^{\kappa T}-|x|e^{\kappa t})}{(1+2D)e^{2\kappa
T}-2De^{2\kappa t}}\cdot\mbox{sign}(x),\\
u_{2}(t,x)  &  =\sqrt{2D}\frac{2\pi}{\sqrt{L\hat{L}}}\frac{2e^{\kappa
t}(e^{\kappa T}-|x|e^{\kappa t})}{(1+2D)e^{2\kappa T}-2De^{2\kappa t}}%
\cdot\mbox{sign}(x),
\end{align*}
where $\mbox{sign}(x)\doteq1$ if $x\geq0$ and $\mbox{sign}(x)\doteq-1$ if
$x<0$.

In the numerical simulation, we set $D=1$, the initial condition
$(t,x)=(0,0.05)$, and the terminal time $T=1$. One can calculate $\hat
{L}=9.83999$ and $\kappa=0.407728$. We used a predictor-corrector Euler scheme
to simulate the trajectories of (\ref{Eq:LangevinEquation2}) and the
associated controlled SDE. For reasons that will be discussed in Appendix
A, a direct numerical approximation scheme was
adopted instead of other techniques such as multiscale integrator or
projective integrator methods (see \cite{VE,ELiuVE,GKK}). By Theorem
\ref{T:DirectScheme} we know that the error in the Euler approximation is
bounded by a term of order $\Delta\epsilon/\delta^{2}$, where $\Delta$ is the
time discretization step. For each choice of $\epsilon$ and $\delta$ we chose
$\Delta$ so that the aforementioned error bound is of the order $0.001$. In
other words, we set $\Delta=0.001\cdot\delta^{2}/\epsilon$. Hence, as
$\epsilon$ gets smaller the discretization step $\Delta$ becomes smaller as
well. Even though by extensive experimentation we found that in general the
choice of $\Delta$ was crucial for obtaining accurate results, for the
periodic example studied here the requirement can be relaxed and a coarser
discretization can still lead to accurate and stable results.

Simulations were done using parallel computing in the C programming language.
We used Mersenne Twister \cite{MatsumotoNishimura1998} for the random number
generator, with a sample size of $N=10^{7}$. The measure for comparing
different schemes is the relative error of the estimators, which is defined
as
\[
\mbox{relative error} \doteq\frac{\mbox{standard deviation of the estimator}}%
{\mbox{expected value of the estimator}}%
\]
The smaller the relative error the more efficient the estimator.
Since in practice both the standard deviation and the expected value
of an estimator are typically unknown, empirical relative error is
often used for measurement. In other words, the expected value of
the estimator will be replaced by the empirical sample mean, and the
standard deviation of the estimator will be replaced by the
empirical sample standard error. In order to distinguish among the
different Monte Carlo procedures, we denote by
$\hat\rho_{i}(\epsilon)$ the empirical relative error of
$\hat\theta_{i}(\epsilon)$ for $i=0,1,2$. We would like to point out
that the expected value in the denominator [which is always
$\theta(\epsilon)$ due to unbiasedness] is replaced by $\hat{\theta
}_{1}(\epsilon)$, regardless of $i$. The reason is that $\hat{\theta}%
_{1}(\epsilon)$ is the most accurate estimate of $\theta(\epsilon)$.

The numerical results are summarized in Table \ref{Table1a}. As suspected, the
estimator $\hat{\theta}_{1}(\epsilon)$ significantly outperforms both the
standard Monte Carlo estimator $\hat{\theta}_{0}(\epsilon)$ and the estimator
$\hat{\theta}_{2}(\epsilon)$ which corresponds to the change of measure purely
based on the homogenized HJB equation. In particular, the estimator
$\hat{\theta}_{1}(\epsilon)$ seems to be of bounded relative error, which is a
stronger notion of efficiency than asymptotic optimality.

\begin{table}[th]
\begin{center}
\begin{small}
\begin{tabular}
[c]{|c|c|c|c|c|c|c|c|c|c|}\hline No. & $\epsilon$ & $\delta$ &
$\epsilon/\delta$ & $\hat{\theta}_{0}(\epsilon)$ &
$\hat{\theta}_{1}(\epsilon)$ & $\hat{\theta}_{2}(\epsilon)$ &
$\hat{\rho }_{0}(\epsilon)$ & $\hat{\rho}_{1}(\epsilon)$ &
$\hat{\rho}_{2}(\epsilon )$\\\hline $1$ & $0.25$ & $0.1$ & $2.5$ &
$2.26e-1$ & $2.25e-1$ & $2.26e-1$ & $3.36e-4$ & $1.76e-3$ &
$6.34e-4$\\\hline $2$ & $0.125$ & $0.04$ & $3.125$ & $3.66e-2$ &
$3.65e-2$ & $3.66e-2$ & $8.40e-4$ & $1.80e-3$ & $1.43e-3$\\\hline
%$3$ & $0.0625$ & $0.015625$ & $4$ & $8.71e-4$ & $8.75e-4$ &
%$8.69e-4$ & $1.06e-2$ & $1.29e-3$ & $4.08e-3$\\\hline
$3$ & $0.063$ & $0.016$ & $3.94$ & $9.34e-4$ & $9.33e-4$ & $9.36e-4$
& $3.29e-3$ & $1.85e-3$ & $4.71e-3$\\\hline $4$ &$0.03125$ & $0.007$
& $4.46$ & $6.93e-7$ & $6.87e-7$ & $6.99e-7$ & $4.47e-2$ & $8.00e-4$
& $3.32e-2$\\\hline $5$ & $0.025$ & $0.004$ & $6.25$ & $1.48e-8$ &
$1.61e-8$ & $1.51e-8$ & $6.85e-2$ & $7.55e-4$ & $3.06e-2$\\\hline
$6$ & $0.02$ & $0.002$ & $10$ & $3.08e-10$ & $1.99e-10$ & $1.51e-10$
& $4.09e-1$ & $3.82e-4$ & $4.97e-2$\\\hline $7$ & $0.015$ & $0.0013$
& $11.54$ & $7.60e-14$ & $1.37e-13$ & $1.07e-13$ & $2.53e-1$ &
$3.01e-4$ & $1.86e-1$\\\hline
\end{tabular}
\end{small}
\end{center}
\caption{Comparison table for periodic case}%
\label{Table1a}%
\end{table}

\subsection{Simulation Results for Random Case}

\label{SS:SimulationResultsRandomCase}

In this section we test the performance of the proposed estimator in the case
of a random environment by estimating an exit probability. In particular, we
again consider the first order Langevin equation (\ref{Eq:LangevinEquation2})
in one dimension and wish to estimate the exit probability
\[
\theta(\epsilon)\doteq\mathrm{P}\left[  X^{\epsilon}(\tau^{\epsilon}%
)=x^{+}|X^{\epsilon}(0)=x\right]  ,
\]
where the exit time $\tau^{\epsilon}$ is defined by
\[
\tau^{\epsilon}\doteq\inf\left\{  t\geq0:X^{\epsilon}(t)\notin(x^{-}%
,x^{+})\right\}  \text{ with }x\in(x^{-},x^{+}).
\]
As in the periodic case, we compare the estimator proposed in Section
\ref{SS:ISRandomMedia} [again denoted by $\hat{\theta}_{1}(\epsilon)]$ with
the standard Monte-Carlo estimator $\hat{\theta}_{0}(\epsilon)$ and with the
estimator $\hat{\theta}_{2}(\epsilon)$ that corresponds to the change of
measure based just on the homogenized HJB equation.

We consider $V(x)=x$ and $Q(y)$ to be a zero mean Gaussian random field with
covariance function
\[
\mathrm{E}^{\nu}\left[  Q(x)Q(y)\right]  =\exp[-\left\vert x-y\right\vert
^{2}].
\]
It follows from Theorem \ref{T:HomogenizationThm1D} that the effective drift
$r(x)$ and diffusivity $q$ are respectively
\[
r=-\kappa,~~~q=2D\kappa,~~~\mbox{where }\kappa=\frac{1}{K\hat{K}}.
\]
In this case the limiting HJB equation is the time independent version of
(\ref{Eq:1DHJBeqn}) defined on the interval $(x^{-},x^{+})$
\[
\kappa U^{\prime}(x)+\kappa D\left\vert U^{\prime}(x)\right\vert ^{2}=0,
\]
with the boundary condition $U(x)=0$ if $x=x^{+}$ and $U(x)=\infty$ if
$x=x^{-}$. We consider the case $D=1$, initial point $(t,x)=(0,0)$ and
$x^{\pm}=\pm0.5$. The maximal viscosity solution (which is also the maximal
classical sense subsolution) to this HJB equation is
\[
U(x)=\frac{1}{D}[x^{+}-x],~~~x\in(x^{-},x^{+}).
\]
While this example is not over a finite time interval, the proof of Theorem
\ref{T:UniformlyLogEfficientRegime1} can be adapted as in \cite{DupuisWang2}
to yield the analogous results. The importance sampling estimator $\hat
{\theta}_{1}(\epsilon)$ corresponds to the control
\[
u_{1}(y)=\frac{\sqrt{2D}}{\hat{K}D}e^{Q(y)/D},
\]
while the estimator $\hat{\theta}_{2}(\epsilon)$ corresponds to the constant control
\[
u_{2}=\sqrt{\frac{2D}{K\hat{K}}}\frac{1}{D}.
\]
Following the notation of the periodic case, we summarize the simulation
results in Table \ref{Table2a}. We used the randomization method to simulate
the Gaussian random environment. See \cite{KramerKurbanmuradovSabelfeld2007}
for an exposition on the simulation of Gaussian random fields.
\begin{table}[th]
\begin{small}
\begin{center}
\begin{tabular}
[c]{|c|c|c|c|c|c|c|c|c|c|}\hline
No. & $\epsilon$ & $\delta$ & $\epsilon/\delta$ & $\hat{\theta}_{0}(\epsilon)$
& $\hat{\theta}_{1}(\epsilon)$ & $\hat{\theta}_{2}(\epsilon)$ & $\hat{\rho
}_{0}(\epsilon)$ & $\hat{\rho}_{1}(\epsilon)$ & $\hat{\rho}_{2}(\epsilon
)$\\\hline
$1$ & $0.25$ & $0.1$ & $2.5$ & $1.38e-1$ & $1.38e-1$ & $1.38e-1$ & $7.88e-4$ &
$1.59e-4$ & $8.09e-4$\\\hline
$2$ & $0.125$ & $0.04$ & $3.125$ & $1.28e-2$ & $1.31e-2$ & $1.28e-2$ &
$2.27e-3$ & $4.91e-3$ & $2.61e-3$\\\hline
$3$ & $0.0625$ & $0.018$ & $3.472$ & $6.02e-4$ & $6.13e-4$ & $5.89e-4$ &
$1.12e-2$ & $5.93e-3$ & $1.32e-2$\\\hline
$4$ & $0.05$ & $0.01$ & $5$ & $2.38e-5$ & $2.30e-5$ & $2.22e-5$ & $6.70e-2$ &
$8.89e-3$ & $9.97e-2$\\\hline
$5$ & $0.04$ & $0.007$ & $5.72$ & $5.5e-6$ & $5.93e-6$ & $4.86e-6$ & $1.25e-1$
& $5.53e-2$ & $1.05e-1$\\\hline
$6$ & $0.025$ & $0.004$ & $6.25$ & $-$ & $7.82e-10$ & $1.26e-09$ & $-$ &
$3.86e-2$ & $5.87e-1$\\\hline
\end{tabular}
\end{center}
\end{small}
\caption{Comparison table for random case with negative drift}%
\label{Table2a}%
\end{table}

As in the periodic case, the estimator $\hat{\theta}_{1}(\epsilon)$
outperforms both $\hat{\theta}_{0}(\epsilon)$ and $\hat{\theta}_{2}(\epsilon
)$. Comparing Tables \ref{Table1a} and \ref{Table2a}, the reader may wonder
why we did not try combinations of $(\epsilon,\delta) $ with larger ratio
$\epsilon/\delta$. Extensive empirical studies showed that the direct
numerical scheme that we chose to simulate from the SDE was much more
sensitive in honoring the rule for choosing the discretization step
$\Delta=0.001\cdot\delta^{2}/\epsilon$ in the random case than it was for the
periodic case. So, due to practical limitation on the computational budget, we
had to limit to the values reported in Table \ref{Table2a} in order to obtain
meaningful results. Note that this significant computational burden required
to produce samples is independent of the scheme, and thus provides further
impetus for the development of the theoretically best algorithms.

We also ran simulations for the model
\begin{equation}
dX^{\epsilon}(t)=\left[  -\frac{\epsilon}{\delta}Q^{\prime}\left(
\frac{X^{\epsilon}(t)}{\delta}\right)  -x\right]  dt+\sqrt{\epsilon}\sqrt
{2D}dW(t),\hspace{0.2cm}X^{\epsilon}(0)=x.\label{Eq:LangevinEquation4}%
\end{equation}
where $Q$ is the same random field as before and the goal is again to estimate
the exit probability
\[
\theta(\epsilon)\doteq\mathrm{P}\left[  X^{\epsilon}(\tau^{\epsilon}%
)=x^{+}|X^{\epsilon}(0)=x\right]  ,
\]
with now $x^{-}=0$, $x^{+}=0.8$ and $x=0.1$.

Notice that this corresponds to the (random) potential function that
is drawn in Figure \ref{F:FigureWithRandomPotential}. The difference
between this and the previous model is that here there exists a
stable point, in particular at $x=0$,  in the domain of attraction
(see Figure \ref{F:FigureWithRandomPotential}). The results were
qualitatively the same as in Table \ref{Table2a} and are presented
in Table \ref{Table3a}.
\begin{table}[th]
\begin{center}%
\begin{small}
\begin{tabular}
[c]{|c|c|c|c|c|c|c|c|c|c|}\hline
No. & $\epsilon$ & $\delta$ & $\epsilon/\delta$ & $\hat{\theta}_{0}(\epsilon)$
& $\hat{\theta}_{1}(\epsilon)$ & $\hat{\theta}_{2}(\epsilon)$ & $\hat{\rho
}_{0}(\epsilon)$ & $\hat{\rho}_{1}(\epsilon)$ & $\hat{\rho}_{2}(\epsilon
)$\\\hline
$1$ & $0.25$ & $0.1$ & $2.5$ & $1.56e-1$ & $1.56e-1$ & $1.56e-1$ & $7.37e-4$ &
$4.93e-4$ & $5.77e-4$\\\hline
$2$ & $0.125$ & $0.04$ & $3.125$ & $2.35e-2$ & $2.39e-2$ & $2.35e-2$ &
$2.00e-3$ & $6.53e-3$ & $1.43e-3$\\\hline
$3$ & $0.0625$ & $0.018$ & $3.472$ & $2.25e-3$ & $2.32e-3$ & $2.25e-3$ &
$6.45e-3$ & $3.66e-2$ & $3.14e-2$\\\hline
$4$ & $0.03125$ & $0.008$ & $3.91$ & $5.03e-5$ & $2.78e-5$ & $4.36e-5$ &
$8.08e-2$ & $8.91e-2$ & $4.54e-1$\\\hline
$5$ & $0.025$ & $0.006$ & $4.17$ & $1.38e-5$ & $5.23e-6$ & $8.91e-6$ &
$2.25e-1$ & $7.79e-2$ & $1.89e-1$\\\hline
$6$ & $0.02$ & $0.0045$ & $4.44$ & $2.0e-7$ & $3.07e-7$ & $3.11e-7$ &
$4.61e-1$ & $1.16e-2$ & $2.01e-1$\\\hline
\end{tabular}
\end{small}
\end{center}
\caption{Comparison table for random case with a rest point}%
\label{Table3a}%
\end{table}

%\appendix

\section*{Appendix A. Numerical Schemes for Multiscale Problems}

\label{SS:NumericalConvergence} Assume for simplicity the periodic setup, and
denote by $Y_{n}^{\epsilon}$ the numerical approximation to $X_{t}^{\epsilon}$
provided by a direct scheme with weak order of convergence $p$. For such a
numerical scheme one has the following error bound, whose proof follows from
standard arguments, e.g. \cite{ELiuVE}, and will not be repeated here.

\begin{theorem}
\label{T:DirectScheme} If $\epsilon,\delta$ and the discretization step
$\Delta$ are such that $\Delta\epsilon/\delta^{2}\ll1$, then for every $T>0$
and every smooth function $f$ with compact support, there exist constants
$C_{0}<\infty$ and $h_{0}>0$ that are independent of $\epsilon,\delta$ and
$\Delta$ such that for all $\Delta,\epsilon,\delta$ satisfying $\Delta
\epsilon/\delta^{2}<h_{0}$,%
\begin{equation}
\sup_{n\leq T/\Delta}|\mathrm{E}_{x_{0}}f(X_{t_{n}}^{\epsilon})-\mathrm{E}%
_{x_{0}}f(Y_{n}^{\epsilon})|\leq C_{0}\hspace{0.1cm}\left(  \Delta
(\epsilon/\delta^{2})\right)  ^{p}. \label{Eq:DirectSchemeBound}%
\end{equation}
%
%\esquare

\end{theorem}

The bound (\ref{Eq:DirectSchemeBound}) illustrates the computational
difficulty in approximating multiscale diffusions. Fix an error tolerance
level $\zeta$. The error bound (\ref{Eq:DirectSchemeBound}) indicates that a
time step of order
\[
\Delta=O\left(  \frac{\delta^{2}}{\epsilon}\zeta^{1/p}\right)
\]
is needed. The corresponding computational cost per unit of time is $1/\Delta
$, which becomes more expensive as $\delta$, $\epsilon$ and $\delta/\epsilon$
become smaller.

Of course, the increasing computational cost highlights the importance of
important sampling or other fast simulation techniques for treating these
kinds of problems. Since our interest in this work is to study importance
sampling and not the numerical methods, we do not elaborate here much on the
approximation aspect of the problem.
Numerical methods such as multiscale integrator methods \cite{VE,ELiuVE}, and
projective integrator methods \cite{GKK, PapavasiliouLevrekidis} have been
proposed to efficiently simulate systems with widely separated time scales.
These methods turn out to be less costly than direct approximation methods,
especially when $\delta$ is small. For example, in the case of equation
(\ref{Eq:LDPandA1}) and for $\delta=\epsilon^{3/2}$, perturbation analysis
\cite{ELiuVE} shows that these methods are less costly than direct
approximation when $\epsilon\ll\zeta$. Observe that in the examples of Section
\ref{S:SimulationResults}, the values of $\epsilon$ and $\delta$ were not
smaller than the overall tolerance error $\zeta$. Hence, we chose to use
direct approximation schemes and to rely on parallel computing to carry out
the simulation.

\section*{Acknowledgments}

We would like to thank the Center for Computation and Visualization (CCV) at
Brown University for making available to us their high performance computing center.

%\bibliography{bibfile}

%%%%%%%%%%%%%%%%%%%%%%%%%%%%%%%%%%%%%%%%%%%%%%%%%%%%%%%%%%%%
%%%%%%%%%%%%%%%%%%%%%%%%%%%%%%%%%%%%%%%%%%%%%%%%%%%%%%%%%%%%%%%%%

\end{document}